\documentclass[a4paper,11pt]{amsart}
\usepackage{amsmath,amssymb,amsthm}
\usepackage{hyperref}
\usepackage{stmaryrd}
\usepackage{xcolor}
\usepackage[T1]{fontenc}
\usepackage[utf8]{inputenc}
\usepackage{xurl}
\usepackage{amscd}
\usepackage{bm}
\usepackage{amsbsy}
\usepackage{wrapfig}
\usepackage{comment}
\let\cal\mathcal

\newcommand{\CE}{\mathrm{CE}}
\newcommand{\fg}{\mathfrak{g}}
\newcommand{\rE}{\mathrm E}

\newcommand{\CC}{\mathbb{C}}

\newcommand{\OO}{\mathcal{O}}
\newcommand{\fd}{\mathfrak{d}}

\newcommand{\Cur}{\mathrm{Cur}}
\newcommand{\kap}{\kappa}
\newcommand{\LCA}{\mathrm{LCA}}
\newcommand{\PVA}{\mathrm{PVA}}

\newcommand{\cE}{\mathcal E}

\newtheorem{dfn}{Definition}[section]
\newtheorem{exa}[dfn]{Example}

\newtheorem{rmk}[dfn]{Remark}
\newtheorem{prop}[dfn]{Proposition}

\newtheorem{thm}{Theorem}[section]
\newcommand{\ch}{\mathrm{ch}}

\newcommand{\Chom}{\mathrm{Chom}}
\newcommand{\id}{\mathrm{id}}
\newcommand{\fsl}{\mathfrak{sl}}
\newcommand{\BBD}{\mathbb{D}}
\newcommand{\fz}{\mathfrak{z}}
\newcommand{\cT}{\mathcal T}
\newcommand{\Span}{\mathrm{Span}}

\newlength{\minitwocolumn}
\newcommand{\beq}{\begin{equation}}
\newcommand{\eeq}{\end{equation}}
\newcommand{\bea}{\begin{eqnarray*}}
\newcommand{\eea}{\end{eqnarray*}}
\newcommand{\beqa}{\begin{eqnarray}}
\newcommand{\eeqa}{\end{eqnarray}}

\newcommand{\End}{\mathrm{End}}
\newcommand{\cV}{\mathcal{V}}
\newcommand{\cW}{\mathcal W}
\newcommand{\cS}{\mathcal S}
\newcommand{\free}{\mathrm{free}}
\newcommand{\wN}{\widetilde{\cN}}

\def\bZ{{\mathbb{Z}}}

\newcommand{\calO}{{\cal O}}

\newcommand{\SL}{\mathrm{SL}}
\newcommand{\CYB}{\mathrm{CYB}}

\newcommand{\ev}{{\rm ev}}

\newcommand{\Vir}{\mathrm{Vir}}

\newcommand{\courant}[2]{{[{{#1},{#2}}]_D}}

\newcommand{\odd}{\mathrm{odd}}

\newcommand{\cL}{\mathcal L}

\newcommand{\bracket}[2]{\langle #1,\,#2\rangle}

\newcommand{\Sym}{\mathrm{Sym}}

\newcommand{\cA}{\mathcal{A}}
\newcommand{\cM}{\mathcal{M}}

\newcommand{\cD}{\mathcal{D}}
\newcommand{\op}{\mathrm{op}}

\newcommand{\cN}{\mathcal N}
\newcommand{\vac}{\mathbf{1}}
\newcommand{\Id}{\mathrm{Id}}

\title{Toward Shifted symplectic Chiralization: Deformation of Courant Algebroids}

\author{Wenda Fang}
\address{Research Institute for Mathematical Sciences, Kyoto University,
 Kyoto 606-8502 JAPAN}
\email{wenda@kurims.kyoto-u.ac.jp}

\author{Noriaki Ikeda}
\address{Department of Mathematical Sciences, Ritsumeikan University,
 Kusatsu 525-8577 JAPAN}
\address{Research Institute for Mathematical Sciences, Kyoto University,
 Kyoto 606-8502 JAPAN}
 \email{nikeda@se.ritsumei.ac.jp}

\thanks{}
\date{\today}

\begin{document}
\begin{abstract}
We study the chiralization of Poisson--Lie theory through graded symplectic geometry. We compare Courant deformations and classical $r$- and $R$-matrix constructions, obtaining integrable hierarchies and free field realizations. We give a variational description of the obstruction to Lie conformal bialgebra compatibility and recover Liberati's obstruction in the Virasoro case. We also introduce a vertex bialgebra analogue of quantized universal enveloping algebras and construct explicit quantizations to all orders of Cartan type Lie conformal bialgebra structures on $\operatorname{Cur}(\mathfrak{sl}_2)$. These results contribute to the first author's broader programme of developing shifted symplectic chiralization as a geometric method for constructing chiral analogues of Poisson-Lie structures via $n$-shifted symplectic geometry, and exploring its connections with AKSZ theory and quantization.
\end{abstract}
\maketitle

\rm

\setcounter{section}{0}

\section{Introduction}
Poisson--Lie theory provides a common language for classical integrable systems and quantum groups. Its central constructions classical $r$-matrices, Lie bialgebras, Manin triples, and Drinfeld doubles connect Poisson geometry,  Hamiltonian dynamics, and quantization, see~\cite{Dr83,Dr87}. Extending these relations to local fields requires algebraic structures that incorporate derivatives into Poisson brackets and Hamiltonian equations. Poisson vertex algebras (PVAs)~\cite{K15,Li04, BDSK09} provide such a setting: their $\lambda$-brackets encode local Poisson brackets, with powers of $\lambda$ recording derivatives of the delta distribution. Lie conformal algebras play the role of Lie algebras in this local setting.

The present paper studies two aspects of this picture for Lie conformal algebras: the realization of $R$-deformed brackets inside doubles of Lie bialgebroids, and the quantization of conformal cobrackets by coproducts on universal enveloping vertex algebras. Our main results are explicit constructions associated with two families of chiralized $r$-matrices on $\Cur(\mathfrak{sl}_2)$. These constructions relate Courant geometry to conformal algebraic structures and yield applications to Hamiltonian systems, free field realizations, and vertex bialgebra quantization.

For background on arc spaces and their relation to vertex algebras, we refer to~\cite{AM}. The geometric approach developed in~\cite{Fang25,Fang26} studies PVAs using infinite dimensional Poisson geometry on arc spaces and shifted symplectic structures. It provides a variational master equation description of local Poisson brackets and a classical $R$-matrix method for constructing compatible PVA structures and integrable systems. In the present paper, we use this approach to compare conformal $r$-matrices, operator $R$-matrices, and Courant $N$-twists, and to examine their consequences for Hamiltonian PDEs, vertex algebras, and vertex bialgebras. We establish explicit relations among these constructions and apply them to examples in which the geometric and algebraic descriptions lead to concrete results.

One motivation comes from quantization. For an ordinary Lie bialgebra $(\fg,\delta)$, the universal enveloping algebra $U(\fg)$ has a primitive cocommutative coproduct, while a quantized enveloping algebra deforms that coproduct so that its first order skew part recovers $\delta$. A Lie conformal algebra $L$ similarly has a universal enveloping vertex algebra $U(L)$ with a canonical cocommutative coproduct. The corresponding conformal quantization problem is to recover a Lie conformal cobracket on $L$ as the first order skew part of a noncocommutative deformation of this coproduct. We give a necessary first order compatibility condition and, for a Cartan type family on $\Cur(\mathfrak{sl}_2)$, an explicit all orders solution that keeps the underlying vertex algebra fixed. The dependence on derivatives distinguishes this construction from the ordinary Lie algebraic case.

We now describe the classical constructions and their applications, followed by the bialgebra compatibility and quantization results.

\subsection*{From conformal $r$-matrices to Hamiltonian systems}
Liberati~\cite{Lib08} introduced Lie conformal bialgebras and conformal classical $r$-matrices. For quadratic current Lie conformal algebras, tensor solutions can be converted into operator solutions by contraction with the invariant conformal bilinear form. This tensor to  operator mechanism is related to~\cite[Corollary~3.3]{HB20}; we give a direct proof in the conventions of~\cite{Fang25} and use it to connect the two deformation procedures.
\begin{thm}\label{thm:CCYBTOOCYB-intro}
Let $(\fg,\kappa)$ be a finite dimensional quadratic Lie algebra over $\CC$, and let $L:=\Cur(\fg)$ carry the induced nondegenerate symmetric invariant conformal bilinear form $\kappa_\lambda$. If $r=\sum_i a_i\otimes b_i\in L\otimes_{\CC}L$ is a skew-symmetric solution of the conformal classical Yang--Baxter equation, define
\[
R_r:L\longrightarrow L,\qquad R_r(x):=\sum_i\kappa_\partial(b_i,x)a_i,
\]
where the substituted $\partial$ acts on $a_i$. Then $R_r$ is a $\CC[\partial]$-module endomorphism and satisfies
\[
[R_r(x)_\lambda R_r(y)]=R_r\bigl([R_r(x)_\lambda y]+[x_\lambda R_r(y)]\bigr)
\]
for all $x,y\in L$. Consequently, $[x_\lambda y]_{R_r}:=[R_r(x)_\lambda y]+[x_\lambda R_r(y)]$ defines a Lie conformal algebra structure on $L$.
\end{thm}

For $L=\Cur(\mathfrak{sl}_2)$, the corresponding $R$-deformed brackets yield concrete Hamiltonian structures on $\mathcal V=\operatorname{Sym}(L)$. In one family, the Hamiltonian operator takes the form $H_{R_c}=H^{(0)}+cH^{(1)}$, so that $H^{(0)}$ and $H^{(1)}$ are compatible. On a localization of $\mathcal V$, we construct the Lenard-Magri sequence $H_n=\int\frac{e(e^{(1)})^n}{(\alpha e-h)^{n+1}}$, $n\geq0$, and obtain an explicit hierarchy of commuting Hamiltonian flows. For a second, Cartan type family, $R_h(e)=R_h(f)=0$ and $R_h(h)=p(\partial)h$ for an odd polynomial $p$, we exhibit a central differential subalgebra after localization and obtain a decomposition of the resulting PVA into a noncentral factor and a factor with zero $\lambda$-bracket. These examples exhibit both Hamiltonian hierarchies and structural decompositions arising from conformal $r$-matrices.

\subsection*{Courant deformations, doubles, and vertex algebra realizations.}

A second route to compatible Poisson structures comes from Courant algebroids. The relation between Courant geometry and local current algebras was established by Alekseev and Strobl~\cite{AS05}, who showed that the current brackets of two dimensional sigma models are governed by the Courant-Dorfman bracket and its natural pairing. This connection was subsequently developed in the QP-manifold framework by Ikeda and Koizumi~\cite{IK13} and Ikeda and Xu~\cite{IX13,IX13v}. In the degree-two shifted symplectic description of a Courant algebroid, the structure is encoded by a degree three Hamiltonian $\Theta$. A weak deforming tensor in the sense of Kosmann-Schwarzbach~\cite{KS11} gives a compatible Hamiltonian $\Theta_{\mathcal N}=\{\widetilde{\mathcal N},\Theta\}$. Over a point, the corresponding affine PVA brackets form a compatible Poisson pencil. We compare this geometric construction with the operator $R$-matrix deformation. For a graded geometric treatment of doubles of Lie bialgebroids, we refer to Voronov~\cite{Vor01}.

For $\mathfrak{sl}_2$, the comparison becomes explicit. The $N$-deformed cotangent double contains a Lie conformal subalgebra $\mathcal L_-$ isomorphic to the $R$-deformed algebra $L_{R_c}$ for $c\neq0$. After removing the central terms, the resulting algebra, together with $\mathcal L_+$ and $\mathcal L_-$, forms a conformal Manin triple. Thus the $R$-deformed algebra is realized as one half of a conformal double arising from a Courant deformation. We refer to this structure as a \emph{chiralized Drinfeld double}; a systematic definition and its relation to $n$-shifted symplectic chiralization will be developed in a forthcoming single authored paper by the first author. The corresponding level one universal affine vertex algebra carries a conformal structure of central charge $6$ and embeds into $\mathcal D^{\mathrm{ch}}_{SL_2,0}$, the vertex algebra of global chiral differential operators on $SL_2$ at level zero.

We also examine the vertex algebras obtained by enveloping the $R$-deformed Lie conformal algebras. For the Cartan family, an additional conformal $2$-cocycle, independent of the pullback of the standard affine cocycle, yields a central extension whose vertex algebra admits a free field factorization $V_{p,k,\kappa}\simeq\mathcal S_\kappa\otimes\mathcal H^{(p)}_{2(k+2)}$, where $\mathcal S_\kappa$ is the rank one Weyl vertex algebra and $\mathcal H^{(p)}_\ell$ is generated by $p(\partial)\phi$ inside the rank-one Heisenberg vertex algebra with $[\phi_\lambda\phi]=\ell\lambda\mathbf1$. When $p(t)=t$, this description also identifies free field factors appearing in the free field limit of the minimal $W$-algebra of $\mathfrak{sl}_3$, see~\cite{ACKL17,CL22}. The same $R$-matrix construction therefore leads both to Hamiltonian structures and to explicit vertex algebra realizations.

\subsection*{Bialgebra compatibility and its geometric obstruction.}The comparison with Lie conformal bialgebras also has an obstruction theoretic aspect. In ordinary Lie bialgebra theory, compatibility of the bracket and cobracket is expressed by a mixed master equation. Using the variational $n$-shifted symplectic formalism of~\cite{Fang26,FangPhD}, we formulate an analogous mixed term in the conformal setting and compare it with Liberati's compatibility obstruction. In the Virasoro case, Theorem~\ref{thm:virob} identifies the two nonzero obstructions, up to the sign convention for the variational bracket. This gives a concrete verification of the geometric description in that example; it does not assert a general equivalence of the two obstruction theories.

\subsection*{From Lie conformal bialgebras to vertex bialgebra quantization.} The final part of the paper returns to the quantum group motivation. Recall that $U(L)$ carries a canonical cocommutative vertex bialgebra structure~\cite{Li07,HLX22}, with primitive elements $P(U(L))=L$; see~\cite{BG22}. The Etingof--Kazhdan theorem establishes quantization for ordinary Lie bialgebras~\cite{EK96,EK98}; here we investigate the analogous coproduct deformation problem for Lie conformal algebras. Our first result in this direction extracts the conformal $1$-cocycle condition from a first order vertex bialgebra deformation.
\begin{thm}\label{thm:QG-intro}
Let $A_\hbar$ be a vertex bialgebra, flat over $\BBD:=\CC[\hbar]/(\hbar^2)$, with a fixed isomorphism of vertex bialgebras $A_\hbar/\hbar A_\hbar\simeq U(L)$. For $a\in U(L)$, choose any lift $\widetilde a\in A_\hbar$ and set
\[
\delta_A(a):=\frac{\Delta_\hbar(\widetilde a)-\Delta_\hbar^{\mathrm{op}}(\widetilde a)}{\hbar}\bmod\hbar.
\]
Then $\delta_A$ is well defined. If $\delta_A(L)\subseteq L\otimes L$, its restriction $\delta:L\to L\otimes L$ is skew-symmetric and satisfies the conformal $1$-cocycle condition
\[
\delta([x_\lambda y])=x_\lambda\delta(y)-y_{-\lambda-D_2}\delta(x),\qquad D_2:=\partial\otimes\id+\id\otimes\partial,
\]
for all $x,y\in L$, where the actions on $L\otimes L$ are the diagonal conformal adjoint actions.
\end{thm}

We also record the second order coassociativity condition for extending such a deformation modulo $\hbar^3$ while keeping the vertex algebra and counit fixed. For the Cartan type conformal cobracket on $\Cur(\mathfrak{sl}_2)$, we then construct an explicit all orders quantization by deforming the coproduct.
\begin{thm}\label{thm:Cartan-quantization-intro}
Let $A:=U(\Cur(\mathfrak{sl}_2))$ be the universal level zero affine vertex algebra. For an odd polynomial $p(t)\in\CC[t]$, set $q_p:=p(\partial)h$ and $K_{\hbar,p}:=\exp(\frac{\hbar}{2}q_p)$, with the exponential interpreted $\hbar$-adically in the commutative differential vertex subalgebra generated by $h$. On $A_\hbar:=A[[\hbar]]$, the formulas
\begin{align*}
\Delta_{\hbar,p}(h)=h\otimes\mathbf1+\mathbf1\otimes h,\ 
\Delta_{\hbar,p}(e)=e\otimes K_{\hbar,p}+K_{\hbar,p}^{-1}\otimes e,\ 
\Delta_{\hbar,p}(f)=f\otimes K_{\hbar,p}^{-1}+K_{\hbar,p}\otimes f
\end{align*}
extend to a coassociative, counital, continuous vertex algebra coproduct $\Delta_{\hbar,p}:A_\hbar\longrightarrow(A\otimes_{\CC}A)[[\hbar]].$
Its semiclassical skew part is the Lie conformal cobracket
\[
\delta_p(h)=0,\qquad\delta_p(e)=e\wedge p(\partial)h,\qquad\delta_p(f)=-f\wedge p(\partial)h.
\]
Thus the underlying vertex algebra stays fixed, while its coproduct is noncocommutative for $p\neq0$.
\end{thm}

The Cartan type quantization in Theorem~\ref{thm:CQ}
also suggests a possible connection with the field-theoretic
approach to quantization through defect fusion.
Motivated by the holomorphic-topological BF theory
considered in~\cite{Zeng26}, we discuss in Remark~9.3
a candidate local interaction functional associated with
the classical cobracket $\delta_p$.
This suggests a possible relationship between our
algebraic quantization and the boundary and defect
fusion constructions studied in
\cite[Sections~5.4 and~6.6]{Zeng26}.
A natural question is whether one can construct a
fusion induced coproduct whose antisymmetric
first order term recovers $\delta_p$.
More generally, one may ask whether the resulting
quantum algebra admits a vertex bialgebra realization
isomorphic to $(A[[\hbar]],\Delta_{\hbar,p})$,
with the isomorphism inducing the prescribed
identification of classical limits.
Establishing such a correspondence would require
compatibility with vertex algebra operations and
translation, together with control of higher order
quantum corrections.
We leave these questions for future investigation.

The results above are concrete applications and extensions of the arc space viewpoint of~\cite{Fang26}. Graded symplectic geometry organizes the comparison with Courant structures and the formulation of the variational obstruction, while the subsequent developments in Hamiltonian theory, conformal algebra, and quantization are carried out in the languages best suited to each. The first author's forthcoming single authored paper will develop a general theory of $n$-shifted symplectic chiralization and chiralized doubles.

This paper is organized as follows. Section~\ref{sec:LCALCBA} recalls Lie conformal algebras, Lie conformal bialgebras, and vertex bialgebras. Section~\ref{sec:QPSCA} recalls graded symplectic structures, Courant algebroids, and the associated current PVA construction. Section~\ref{sec:NDCA} studies Courant $N$-deformations and compatible PVA brackets in the point base setting. Section~\ref{sec:CRM} relates conformal $r$-matrices to operator $R$-matrices and compares the resulting Lie conformal algebras with the deformed double of Lie bialgebroids. Sections~\ref{sec:APPIS} and~\ref{sec:APVAS} study integrable systems, central extensions, and free field realizations. Section~\ref{sec:OQVOS} compares the Virasoro compatibility obstruction with its variational shifted symplectic description. Section~\ref{Sec:QT} studies coproduct deformations and the all orders Cartan quantization.
\section*{Acknowledgments}
The first author thanks Tomoyuki Arakawa for discussions on this work and related earlier work at SLMath, SNU, and Luminy, for his encouragement, and for the opportunity to present his work at OIST. The authors also thank Ping Xu, Kai Behrend, and Mathieu Stiénon for discussions during the GAP conference and for their encouragement. The first author is especially grateful to Takahiro Shiota for almost weekly discussions, for listening to his explanations of the $n$-shifted symplectic chiralization formalism, and for his continued encouragement. The first author thanks the organizers of GAP for the opportunity to present his work on various topics. The authors thank Hsuan-Yi Liao and Kaichuan Qi for helpful discussions. The first author also thanks Anne Moreau for discussions at Luminy and Simon Lentner for discussions at RIMS and OIST. The authors are grateful to all the speakers and audience members of their seminar series for their participation and stimulating discussions. Part of this work was completed during Poisson 2026 and GAP 2026. The authors thank the organizers of both conferences for their hospitality and support. This work was supported by JSPS KAKENHI Grant Number 22K03323. The first author used ChatGPT (OpenAI) to assist
with English language editing and proofreading.

\section{Preliminary I: Lie conformal algebra, Lie conformal bialgebras}\label{sec:LCALCBA}
Throughout this paper, we follow the standard conventions for Lie conformal algebras and vertex algebras in \cite{Kac98}. For the theory of classical $r$-matrices on Lie conformal algebras, we refer to \cite{Fang25}. We briefly recall the conformal analogues of Lie coalgebras and Lie bialgebras, following Liberati~\cite{Lib08}. In a forthcoming single authored paper, the first author will refer to Liberati's approach as \emph{Kac-Liberati chiralization}.
\subsection{Conformal linear algebra}
In this section, we recall the notion of a module over a Lie conformal algebra.
\begin{dfn}\label{dfn:LCAMOD}
Let $L$ be a Lie conformal algebra. An $L$-module is a $\CC[\partial]$-module $M$, with $\partial$-action denoted by $\partial_M$, together with a $\CC$-bilinear map $L\otimes M\rightarrow M[\lambda]$, $a\otimes m\mapsto a_\lambda m$, such that, for all $a,b\in L$ and $m\in M$,
\[
(\partial a)_\lambda m=-\lambda a_\lambda m,\qquad a_\lambda(\partial_Mm)=(\partial_M+\lambda)a_\lambda m,\qquad a_\lambda(b_\mu m)-b_\mu(a_\lambda m)=[a_\lambda b]_{\lambda+\mu}m.
\]
\end{dfn}
\begin{exa}\label{exa:LCAMOD}
Let $U$ and $V$ be $L$-modules. Then $M=U\otimes_{\CC}V$ is an $L$-module with
\[
\partial_M=\partial\otimes 1+1\otimes\partial,\qquad a_\lambda(u\otimes v)=(a_\lambda u)\otimes v+u\otimes(a_\lambda v),
\]
for all $a\in L$ and $u\otimes v\in U\otimes V.$
\end{exa}
\begin{dfn}\label{dfn:CLM}
Let $U$ and $V$ be $\CC[\partial]$-modules. A conformal linear map from $U$ to $V$ is a $\CC$-linear map $a:U\rightarrow V[\lambda]$, $u\mapsto a_\lambda(u)$, such that $a_\lambda(\partial^Uu)=(\partial^V+\lambda)a_\lambda(u),$
for all $u\in U$, or equivalently, $\partial^Va_\lambda-a_\lambda\partial^U=-\lambda a_\lambda$. The space of all such maps is denoted by $\Chom(U,V)$ and is a $\CC[\partial]$-module with $(\partial a)_\lambda=-\lambda a_\lambda$. The conformal dual of $U$ is defined by $U^*_{\ch}:=\Chom(U,\CC)$, where $\CC$ is the trivial $\CC[\partial]$-module; equivalently,
\[
U^*_{\ch}=\left\{a:U\rightarrow\CC[\lambda]\mid a\text{ is $\CC$-linear and }a_\lambda(\partial^Uu)=\lambda a_\lambda(u)\text{ for all }u\in U\right\}.
\]
\end{dfn}
\begin{dfn}\label{dfn:CBF}
Let $V$ be a $\CC[\partial]$-module. A conformal bilinear form on $V$ is a $\CC$-bilinear map
\[
(\,\cdot,\cdot\,)_\lambda:V\times V\longrightarrow\CC[\lambda]
\]
such that $(\partial v,w)_\lambda=-\lambda(v,w)_\lambda$ and $(v,\partial w)_\lambda=\lambda(v,w)_\lambda$ for all $v,w\in V$. The conformal bilinear form is called symmetric if $(v,w)_\lambda=(w,v)_{-\lambda}$ for all $v,w\in V$.

Let $L$ be a Lie conformal algebra. A conformal bilinear form on $L$ is called invariant if
\[
([a_\mu b],c)_\lambda=(a,[b_{\lambda-\partial}c])_\mu=-(a,[c_{-\lambda}b])_\mu
\]
for all $a,b,c\in L$.

Given a conformal bilinear form on a $\CC[\partial]$-module $V$, it induces a homomorphism of $\CC[\partial]$-modules
\[
\Phi:V\longrightarrow V^*_{\ch},\qquad v\longmapsto\Phi_v,
\]
defined by $(\Phi_v)_\lambda(w):=(v,w)_\lambda$ for all $v,w\in V$. Following~\cite{Lib08}, we say that an invariant conformal bilinear form $(\,\cdot,\cdot\,)_\lambda$ is non-degenerate if the induced map $\Phi:V\longrightarrow V^*_{\ch}$ is an isomorphism of $\CC[\partial]$-modules.
\end{dfn}
\begin{exa}\label{exa:CBL}
Let $(\fg,(\cdot,\cdot))$ be a finite dimensional quadratic Lie algebra, i.e.,
$(\cdot,\cdot)$ is a non-degenerate symmetric bilinear form on $\fg$ satisfying $([a,b],c)=(a,[b,c])$
for all $a,b,c\in\fg$. Consider the current Lie conformal algebra $\Cur(\fg)=\CC[\partial]\otimes \fg,
\ 
[a_\lambda b]=[a,b],
\  a,b\in\fg,$
extended to all of $\Cur(\fg)$ by sesquilinearity.

Define a $\CC$-bilinear map $(\cdot,\cdot)_\lambda:
\Cur(\fg)\times \Cur(\fg)\longrightarrow \CC[\lambda]$
by
\[
\bigl(p(\partial)a,q(\partial)b\bigr)_\lambda
:=
p(-\lambda)q(\lambda)(a,b),
\qquad
p(\partial),q(\partial)\in\CC[\partial],\ a,b\in\fg.
\]

Then $(\cdot,\cdot)_\lambda$ is a conformal bilinear form. Since
$(\cdot,\cdot)$ is symmetric, $(\cdot,\cdot)_\lambda$ is symmetric.

By sesquilinearity and the invariance of $(\cdot,\cdot)$, it is also invariant.

Finally, let $\{e_i\}_{i=1}^N$ be a basis of $\fg$. Then the conformal Gram
matrix of $(\cdot,\cdot)_\lambda$ with respect to the $\CC[\partial]$-basis
$\{e_i\}_{i=1}^N$ of $\Cur(\fg)$ is the constant matrix $G=\bigl((e_i,e_j)\bigr)_{i,j}$.
Since $(\cdot,\cdot)$ is non-degenerate on $\fg$, we have
$\det G\in\CC^\times$. Therefore the induced map $\Cur(\fg)\longrightarrow \Cur(\fg)^*_{\ch}$
is an isomorphism of $\CC[\partial]$-modules; see~\cite[Section 3.3]{Lib08}. Thus
$(\cdot,\cdot)_\lambda$ is non-degenerate in the sense of
Definition~\ref{dfn:CBF}.
\end{exa}

\subsection{Cohomology of Lie conformal algebra}
As in Poisson--Lie theory, classical $r$-matrices are naturally related to the Chevalley--Eilenberg complex. Kac introduced Lie conformal algebras as a chiral analogue of Lie algebras, and Liberati~\cite{Lib08} developed a chiral analogue of Lie bialgebras using LCA cohomology. This cohomology theory was developed in~\cite{BKV99,DSK09}; we recall the necessary material, together with the basic theory of central extensions of Lie conformal algebras needed in Section~\ref{sec:APVAS}.
\begin{dfn}\label{dfn:COHOMLCA}
Let \(L\) be a Lie conformal algebra and \(M\) an \(L\)-module. For \(n\in\mathbb Z_{\geq0}\), a basic \(n\)-cochain of \(L\) with coefficients in \(M\) is a \(\mathbb C\)-linear map
\[
\gamma:L^{\otimes n}\longrightarrow M[\lambda_1,\ldots,\lambda_n],\qquad a_1\otimes\cdots\otimes a_n\longmapsto\gamma_{\lambda_1,\ldots,\lambda_n}(a_1,\ldots,a_n),
\]
such that, for every \(1\leq i\leq n\) and \(\sigma\in S_n\),
\[
\gamma_{\lambda_1,\ldots,\lambda_n}(a_1,\ldots,\partial a_i,\ldots,a_n)=-\lambda_i\gamma_{\lambda_1,\ldots,\lambda_n}(a_1,\ldots,a_n),
\]
\[
\gamma_{\lambda_{\sigma(1)},\ldots,\lambda_{\sigma(n)}}(a_{\sigma(1)},\ldots,a_{\sigma(n)})=\operatorname{sgn}(\sigma)\gamma_{\lambda_1,\ldots,\lambda_n}(a_1,\ldots,a_n).
\]
The space of basic \(n\)-cochains is denoted by \(\widetilde C^n(L,M)\). We set \(L^{\otimes0}=\mathbb C\), so that \(\widetilde C^0(L,M)\cong M\).
\end{dfn}

Define \(d_{\CE}^{\ch}:\widetilde C^n(L,M)\to\widetilde C^{n+1}(L,M)\) by
\begin{align*}
&(d_{\CE}^{\ch}\gamma)_{\lambda_1,\ldots,\lambda_{n+1}}(a_1,\ldots,a_{n+1})\\
&=\sum_{i=1}^{n+1}(-1)^{i+1}a_{i\,\lambda_i}
\gamma_{\lambda_1,\ldots,\widehat{\lambda_i},\ldots,\lambda_{n+1}}
(a_1,\ldots,\widehat{a_i},\ldots,a_{n+1})\\
&\quad+\sum_{1\leq i<j\leq n+1}(-1)^{i+j}
\gamma_{\lambda_i+\lambda_j,\lambda_1,\ldots,\widehat{\lambda_i},\ldots,\widehat{\lambda_j},\ldots,\lambda_{n+1}}
([a_i{}_{\lambda_i}a_j],a_1,\ldots,\widehat{a_i},\ldots,\widehat{a_j},\ldots,a_{n+1}),
\end{align*}
where a hat denotes omission and \(\gamma\) is extended linearly to polynomial-valued arguments.

Since $(d_{\mathrm{CE}}^{\mathrm{ch}})^2=0$, the cochains
form the basic complex $\widetilde C^\bullet(L,M)
=\bigoplus_{n\in\bZ_{\geq0}}\widetilde C^n(L,M).$
We equip $\widetilde C^\bullet(L,M)$ with a
$\CC[\partial]$-module structure by
\[
(\partial\gamma)_{\lambda_1,\ldots,\lambda_n}
(a_1,\ldots,a_n)
=
\left(\partial_M+\sum_{i=1}^n\lambda_i\right)
\gamma_{\lambda_1,\ldots,\lambda_n}
(a_1,\ldots,a_n),
\]
where $\partial_M$ denotes the action of $\partial$ on $M$.
Since $d_{\mathrm{CE}}^{\mathrm{ch}}\partial
=\partial d_{\mathrm{CE}}^{\mathrm{ch}}$,
the reduced complex is defined by
\[
C^\bullet(L,M)
:=
\frac{\widetilde C^\bullet(L,M)}
{\partial\widetilde C^\bullet(L,M)}
=
\bigoplus_{n\in\bZ_{\geq0}}C^n(L,M),
\qquad
C^n(L,M):=
\frac{\widetilde C^n(L,M)}
{\partial\widetilde C^n(L,M)}.
\]
For simplicity, we henceforth write
$d$ for $d_{\mathrm{CE}}^{\mathrm{ch}}$.

\begin{rmk}
We compare the Lie conformal cohomology considered above with variational PVA cohomology and the coisson cohomology studied by Butson-Nair~\cite{BN26}. Let $L$ be a Lie conformal algebra which is free of finite rank over $\CC[\partial]$, and let $\cV:=S^\bullet(L)$ be its universal Poisson vertex algebra. Denote by $C^\bullet_{\LCA}(L,M)$ the reduced Lie conformal cochain complex with coefficients in an $L$-module $M$.

Equip $L\otimes L$ with the diagonal adjoint action and the total translation operator. Since the tensor flip commutes with both actions, we have a decomposition of $L$-modules $L\otimes L=\wedge^2L\oplus S^2L$. Consequently,
\[
H^k_{\LCA}(L,L\otimes L)\cong H^k_{\LCA}(L,\wedge^2L)\oplus H^k_{\LCA}(L,S^2L),\qquad k\geq 0.
\]
Moreover, since $S^2L$ is a direct summand of $\cV$ as an $L$-module, the comparison theorem of ~\cite[Theorem~3.13(a)]{BDSK20} gives a natural split injection $H^k_{\LCA}(L,S^2L)\hookrightarrow H^k_{\LCA}(L,\cV)\cong H^k_{\PVA,\mathrm{var}}(\cV,\cV).$
On the other hand, ~\cite[Corollary~5.45]{BN26} prove that the inclusion of the variational Poisson deformation complex into the coisson deformation complex is a quasi-isomorphism when the underlying commutative differential algebra is very smooth. This condition holds for the differential polynomial algebra $\cV$. Taking into account their deformation theoretic degree convention, we obtain $H^k_{\PVA,\mathrm{var}}(\cV,\cV)\cong H^{k-1}_{c}(\cV),\ k\geq 2.$
Thus the symmetric summand $H^k_{\LCA}(L,S^2L)$ admits a natural comparison with coisson cohomology. However, the projection $L\otimes L\to S^2L$ annihilates $\wedge^2L$, and hence this comparison does not detect the skew-symmetric cobrackets arising from conformal $r$-matrices.
\end{rmk}

\begin{exa}\label{exa:LLr} Let $M=L\otimes L$ be equipped with the diagonal conformal adjoint action $a_\lambda(u\otimes v)=[a_\lambda u]\otimes v+u\otimes[a_\lambda v]$ for all $a,u,v\in L$. By Example~\ref{exa:LCAMOD}, $M$ is an $L$-module with $\partial_M=\partial\otimes1+1\otimes\partial=:D_2$, and $\widetilde C^0(L,M)=M$. For $r\in M$, Definition~\ref{dfn:COHOMLCA} gives $(d_{\CE}^{\ch}r)_\lambda(a)=a_\lambda r$. The reduced complex is $C^\bullet(L,M):=\widetilde C^\bullet(L,M)/\partial\widetilde C^\bullet(L,M)$. Since $(\partial\gamma)_\lambda=(D_2+\lambda)\gamma_\lambda$ for $\gamma\in\widetilde C^1(L,M)$, reduction amounts to setting $\lambda=-D_2$. Hence the reduced class of $d_{\CE}^{\ch}r$ is represented by $\delta_r^{\ch}:L\longrightarrow L\otimes L$, where $\delta_r^{\ch}(a):=\left.a_\lambda r\right|_{\lambda=-D_2}$. To make this notation explicit, suppose $r=u\otimes v$. Then $a_\lambda r=\sum_{n\geq0}\frac{\lambda^n}{n!}\bigl((a_{(n)}u)\otimes v+u\otimes(a_{(n)}v)\bigr)$ and hence $\delta_r^{\ch}(a)=\sum_{n\geq0}\frac{(-D_2)^n}{n!}\bigl((a_{(n)}u)\otimes v+u\otimes(a_{(n)}v)\bigr)$, where $D_2(x\otimes y)=\partial x\otimes y+x\otimes\partial y$, so $(-D_2)^n$ acts on the entire tensor in parentheses. The formula extends linearly to every $r\in L\otimes L$. We write simply $\delta_r$ when no confusion can arise. \end{exa}
We will interpret this element $r\in L\otimes L$ as a chiral analogue of a classical $r$-matrix in Section~\ref{sec:CRM}.

Another important example of LCA cohomology is as $M=\CC K.$
\begin{exa}\label{exa:LCAEXT}
Assume that $L$ is free as a $\CC[\partial]$-module. Let $M=\CC K$ be the trivial $L$-module, namely, $\partial_MK=0$ and $a_\lambda K=0$ for all $a\in L$. A basic $2$-cochain $\widetilde\gamma\in\widetilde C^2(L,M)$ is a $\CC$-linear map
 $\widetilde\gamma:L\otimes L\longrightarrow M[\lambda_1,\lambda_2].$
Writing $\widetilde\gamma_{\lambda_1,\lambda_2}(a,b)=\widetilde\omega_{\lambda_1,\lambda_2}(a,b)K$, we have $(\partial\widetilde\gamma)_{\lambda_1,\lambda_2}=(\lambda_1+\lambda_2)\widetilde\gamma_{\lambda_1,\lambda_2}.$
Hence its class in $C^2(L,M)$ is represented by $\omega_\lambda(a,b):=\widetilde\omega_{\lambda,-\lambda}(a,b)\in\CC[\lambda]$. If $\omega$ is a $2$-cocycle, it defines the central extension $\widehat L=L\oplus\CC K$ with $[a_\lambda b]_{\widehat L}=[a_\lambda b]+\omega_\lambda(a,b)K,\  [K_\lambda\widehat L]=0.$
Equivalence classes of such extensions are classified by $H^2(L,M)$~\cite[Theorem~3.1(4)]{BKV99}.
\end{exa}
In Section~\ref{sec:OQVOS}, we briefly explain how $n$-shifted symplectic chiralization can be used to understand Lie conformal bialgebras (LCBAs); here, we recall the definition and give an example.
\begin{dfn}\label{dfn:LCCA}
 A Lie conformal coalgebra is a
$\mathbb C[\partial]$-module $L$ equipped with a $\mathbb C[\partial]$-linear
map $\delta:L\longrightarrow \wedge^2 L$
satisfying the coJacobi identity $(I\otimes\delta)\delta-\tau_{12}(I\otimes\delta)\delta
=(\delta\otimes I)\delta .$
Here $\partial$ acts diagonally on tensor products. A Lie conformal bialgebra
is a triple $(L,[{-}_\lambda{-}],\delta)$ such that $(L,[{-}_\lambda{-}])$ is a
Lie conformal algebra, $(L,\delta)$ is a Lie conformal coalgebra, and $a_\lambda\delta(b)-b_{-\lambda-D_2}\delta(a)
=
\delta([a_\lambda b]),$
for all $a,b\in L$, where the action on $L\otimes L$ is the diagonal conformal
adjoint action. The last identity is the conformal $1$-cocycle condition.
When $L$ is free as a $\CC [\partial]$-module, it is equivalent to saying that $\delta$ corresponds to a $1$-cocycle in the reduced complex of Definition~\ref{dfn:COHOMLCA}.
\end{dfn}
We give an example here. Gubarev and Kozlov~\cite{GK23}
classified the solutions of the conformal classical
Yang-Baxter equation and its weak version on $\Cur(\fsl_2)$.
The conformal classical Yang-Baxter equation may be viewed
as a chiral analogue of the classical Yang--Baxter equation.

\begin{exa}\label{exa:LCBAsl2}
Consider the current Lie conformal algebra \(L=\Cur(\fsl_2)\). As a \(\CC[\partial]\)-module, it has basis \(f,h,e\), and its \(\lambda\)-bracket is induced by the Lie bracket of \(\fsl_2\). Write \(e^{(1)}:=\partial e\) and define $\delta:L\longrightarrow\wedge^2L$
on the basis by \(\delta(e)=0,\ \delta(h)=2e^{(1)}\wedge e,\ \delta(f)=h\wedge e^{(1)}\), extending it \(\CC[\partial]\)-linearly, so that \(\delta(\partial x)=D_2\delta(x)\). Then \(\bigl(\Cur(\fsl_2),[\,\cdot_\lambda\cdot\,],\delta\bigr)\) is a Lie conformal bialgebra. Indeed, a direct computation on the generators \(e,h,f\) shows that \(\delta\) satisfies the coJacobi identity and the compatibility condition $a_\lambda\delta(b)-b_{-\lambda-D_2}\delta(a)=\delta([a_\lambda b]),\ a,b\in L.$
Moreover, this is the coboundary LCBA associated with the chiral classical \(r\)-matrix \(r=\frac{1}{2}e^{(1)}\wedge e\); see Liberati~\cite[Definition~3.3]{Lib08} and Definition~\ref{dfn:CCRM}. A geometric interpretation of this construction will be given in future work of the first author.
\end{exa}

\begin{dfn}\label{dfn:FR}
Let $(L,\partial,[\,\cdot_\lambda\cdot\,])$ be an LCA. An element
$R\in\End_{\CC[\partial]}L$ is called a classical $R$-matrix of $L$ if
and only if the bracket $[X_\lambda Y]_R
=
[R(X)_\lambda Y]+[X_\lambda R(Y)]\ for\ all\ X,Y\in L,$
defines an LCA structure on $L$.
\end{dfn}

\begin{exa}\label{exa:FR}
Let $L=v_k(\fg)$ be the affine Lie conformal algebra at level $k$, and suppose that $\fg=\fg_+\oplus\fg_-$, where $\fg$ is a quadratic Lie algebra equipped with a nondegenerate ad-invariant bilinear form $(\,\cdot,\cdot\,)$, and $\fg_-$ is an isotropic Lie subalgebra with respect to this bilinear form and $\fg_+$ is a Lie subalgebra. Then the factorization $v_k(\fg)=v_k(\fg_+)\oplus v_k'(\fg_-)$ yields a factorizable $R$-matrix. Here $v_k'(\fg_-)$ denotes the current Lie conformal algebra of $\fg_-$, and the central summand is included only in $v_k(\fg_+)$.
\end{exa}

\subsection{Vertex bialgebras}
In this section, we recall the basic definitions of vertex algebras and
vertex bialgebras that will be used in Section~\ref{Sec:QT}. Our main
references are~\cite{HLX22,DSK06}.

\begin{dfn}[{\cite[Definition~1.13]{DSK06}}]\label{dfn:VA}
A vertex algebra is a quintuple $(V,\vac,T,[\,\cdot_\lambda\cdot\,],:\!\cdot\,\cdot\!:)$ satisfying the following conditions:
\begin{itemize}
\item[1.] $(V,T,[\,\cdot_\lambda\cdot\,])$ is a Lie conformal algebra.
\item[2.] $(V,\vac,T,:\!\cdot\,\cdot\!:)$ is a unital differential algebra satisfying quasi-commutativity and quasi-associativity:
\[
:\!ab\!:-:\!ba\!:=\int_{-T}^{0}[a_\lambda b]\,d\lambda,
\]
\[
:\!(:\!ab\!:)c\!:-:\!a(:\!bc\!:)\!:=:\!\left(\int_0^T d\lambda\,a\right)[b_\lambda c]\!:+:\!\left(\int_0^T d\lambda\,b\right)[a_\lambda c]\!:.
\]
\item[3.] The $\lambda$-bracket and normally ordered product satisfy the non-commutative Wick formula:
\[
[a_\lambda:\!bc\!:]=:\![a_\lambda b]c\!:+:\!b[a_\lambda c]\!:+\int_0^\lambda [[a_\lambda b]_\mu c]\,d\mu.
\]
\end{itemize}
\end{dfn}
\begin{dfn}\label{dfn:VBA}
A vertex bialgebra is a vertex algebra $(V,\mathbf{1},\partial,[\,\cdot_\lambda\cdot\,],:\!\cdot\,\cdot\!:)$ equipped with vertex algebra homomorphisms $\Delta:V\longrightarrow V\otimes V$ and $\epsilon:V\longrightarrow\CC$ satisfying $(\Delta\otimes\Id)\Delta=(\Id\otimes\Delta)\Delta,\ (\epsilon\otimes\Id)\Delta=\Id=(\Id\otimes\epsilon)\Delta,$
where $V\otimes V$ carries the tensor product vertex algebra structure and $\CC$ the trivial vertex algebra structure. Explicitly, for all $a,b\in V$,
\[
\Delta(\mathbf{1})=\mathbf{1}\otimes\mathbf{1},\qquad\Delta(\partial a)=(\partial\otimes\Id+\Id\otimes\partial)\Delta(a),
\]
\[
\Delta(:\!ab\!:)=:\!\Delta(a)\Delta(b)\!:,\qquad\Delta([a_\lambda b])=[\Delta(a)_\lambda\Delta(b)]_{V\otimes V},
\]
and
\[
\epsilon(\mathbf{1})=1,\qquad\epsilon(\partial a)=0,\qquad\epsilon(:\!ab\!:)=\epsilon(a)\epsilon(b),\qquad\epsilon([a_\lambda b])=0.
\]
\end{dfn}

\begin{rmk}\label{rmk:TVA}
For the tensor product vertex algebra $V\otimes V$,
\[
Y_{V\otimes V}(u\otimes v,z)=Y_V(u,z)\otimes Y_V(v,z)=\sum_{n\in\mathbb Z}(u\otimes v)_{(n)}z^{-n-1}.
\]
Hence
$(u\otimes v)_{(n)}(\widetilde u\otimes\widetilde v)=\sum_{m\in\mathbb Z}(u_{(m)}\widetilde u)\otimes(v_{(n-m-1)}\widetilde v),$
and therefore
\[
[(u\otimes v)_\lambda(\widetilde u\otimes\widetilde v)]_{V\otimes V}=\sum_{\substack{n\geq0\\m\in\mathbb Z}}\frac{\lambda^n}{n!}(u_{(m)}\widetilde u)\otimes(v_{(n-m-1)}\widetilde v).
\]
In particular, this is generally different from $[u_\lambda\widetilde u]_V\otimes[v_\lambda\widetilde v]_V$.
\end{rmk}
As in the Lie algebra case, every Lie conformal algebra $L$ admits a
universal enveloping vertex algebra, denoted by $U(L)$; see
\cite{Kac98,Primc99}. The following proposition is due to~\cite{Li07}.
\begin{prop}\label{prop:canonical-vertex-bialgebra}
Let $L$ be a Lie conformal algebra over $\CC$. There are unique vertex algebra homomorphisms $\Delta_0:U(L)\longrightarrow U(L)\otimes U(L),
\
\varepsilon_0:U(L)\longrightarrow\CC,$
determined by
$
\Delta_0(x)=x\otimes\mathbf{1}+\mathbf{1}\otimes x,
\
\varepsilon_0(x)=0,
\ x\in L.
$
These maps equip $U(L)$ with a cocommutative vertex bialgebra structure.
\end{prop}

\section{Preliminaries II: shifted symplectic structures and Courant algebroids}
\label{sec:QPSCA}

In this section we recall Roytenberg's description of Courant
algebroids in terms of $2$-shifted symplectic dg varieties and
then explain the current construction leading to the Poisson vertex
algebras used below. Courant algebroids were introduced by Liu,
Weinstein, and Xu~\cite{LWX97}; their degree-$2$ symplectic
description is due to Roytenberg~\cite{Roy01}. We use the Dorfman
bracket convention throughout.

\begin{dfn}
Let $M$ be a smooth complex algebraic variety. A \emph{nonnegatively graded algebraic variety} with body $M$ is a graded locally ringed space $\mathcal M=(M,\mathcal O_{\mathcal M})$ such that $\mathcal O_{\mathcal M}=\bigoplus_{r\geq0}\mathcal O_{\mathcal M}^r$, $\mathcal O_{\mathcal M}^0=\mathcal O_M$, and Zariski locally $\mathcal O_{\mathcal M}|_U\simeq\mathcal O_U\otimes_{\CC}S_{\mathrm{gr}}(W),
$
where $W$ is a finite-dimensional positively graded vector space and $S_{\mathrm{gr}}(W)$ is its free graded commutative algebra. Thus $fg=(-1)^{|f||g|}gf$ for homogeneous sections, with parity given by degree modulo $2$. If the local generators have degrees at most $n$, then $\mathcal M$ has \emph{degree at most $n$}. We write $\mathcal O(\mathcal M):=\Gamma(M,\mathcal O_{\mathcal M})$.
\end{dfn}

\begin{dfn}
The graded tangent sheaf of $\mathcal M$ is $\mathcal T_{\mathcal M}:=\underline{\operatorname{Der}}_{\CC}(\mathcal O_{\mathcal M},\mathcal O_{\mathcal M})$. A \emph{differential} on $\mathcal M$ is a degree-$1$ derivation $Q\in\Gamma(M,\mathcal T_{\mathcal M})^1$ satisfying $Q^2=0$. Thus, for homogeneous local sections $f,g$, $Q(fg)=Q(f)g+(-1)^{|f|}fQ(g),\ [Q,Q]=2Q^2=0.$
The pair $(\mathcal M,Q)$ is a
\emph{nonnegatively graded dg algebraic variety}.
\end{dfn}

All internal gradings are cohomological, and repeated indices are summed. Let $\Omega^1_{\mathcal M}$ be the sheaf of graded K\"ahler differentials; its de Rham differential $d_{\mathrm{dR}}$ has form degree $1$ and internal degree $0$.

\begin{dfn}\label{dfn:QP2}
Let $n\geq1$ and let $\mathcal M=(M,\mathcal O_{\mathcal M})$ be a nonnegatively graded algebraic variety as in Definition~3.1. A symplectic form of degree $n$ on $\mathcal M$ is a form $\omega\in\Gamma(M,\Omega^2_{\mathcal M})^n$ satisfying $d_{\mathrm{dR}}\omega=0$ such that contraction induces an isomorphism of graded sheaves $\omega^\flat:\mathcal T_{\mathcal M}\xrightarrow{\sim}\Omega^1_{\mathcal M}[n],
\  X\longmapsto\iota_X\omega,$
where $(\Omega^1_{\mathcal M}[n])^r=(\Omega^1_{\mathcal M})^{r+n}$. An $n$-shifted symplectic dg algebraic variety is a triple $(\mathcal M,\omega,Q)$, where $\omega$ is such a symplectic form and $Q\in\Gamma(M,\mathcal T_{\mathcal M})^1$ is a derivation satisfying $Q^2=0,\  \mathcal L_Q\omega=0.$
The induced Poisson bracket has degree $-n$: $\{-,-\}_{\mathcal M}:
\mathcal O_{\mathcal M}^p(U)\times\mathcal O_{\mathcal M}^q(U)
\longrightarrow\mathcal O_{\mathcal M}^{p+q-n}(U)$
for every open subset $U\subseteq M$ and all $p,q\geq0$.
\end{dfn}

Since $\mathcal L_Q\omega=0$ and $n+1>0$, the Euler vector field implies that $Q$ has a unique Hamiltonian $\Theta\in\mathcal O(\mathcal M)^{n+1}$. We use the convention $Q=\{\Theta,-\}_{\mathcal M}$. With this convention, $Q^2=0$ is equivalent to the classical master equation $\{\Theta,\Theta\}_{\mathcal M}=0$.
\begin{dfn}\label{dfn:CA}
Let $M$ be a smooth complex algebraic variety and $\cT_M:=\underline{\operatorname{Der}}_{\CC}(\OO_M,\OO_M)$. A \emph{Courant algebroid} over $M$ consists of a finite-rank locally free $\OO_M$-module $\cE$, a symmetric $\OO_M$-bilinear pairing $\bracket{-}{-}:\cE\otimes_{\OO_M}\cE\to\OO_M$ inducing an isomorphism $\cE\xrightarrow{\sim}\cE^*$, an $\OO_M$-linear anchor $\rho:\cE\to\cT_M$, and a $\CC$-bilinear Dorfman bracket $\courant{-}{-}:\cE(U)\times\cE(U)\to\cE(U)$ compatible with restrictions. For every open $U\subseteq M$, all $e,e_1,e_2,e_3\in\cE(U)$ and $f\in\OO_M(U)$, these data satisfy
\begingroup
\setlength{\jot}{2pt}
\begin{align*}
\courant{e_1}{\courant{e_2}{e_3}}
&=\courant{\courant{e_1}{e_2}}{e_3}
+\courant{e_2}{\courant{e_1}{e_3}},\\
\rho(\courant{e_1}{e_2})
&=[\rho(e_1),\rho(e_2)],
\qquad
\courant{e}{e}=\tfrac12\mathcal D\bracket{e}{e},\\
\courant{e_1}{fe_2}
&=f\courant{e_1}{e_2}+(\rho(e_1)f)e_2,\\
\rho(e_1)\bracket{e_2}{e_3}
&=\bracket{\courant{e_1}{e_2}}{e_3}
+\bracket{e_2}{\courant{e_1}{e_3}}.
\end{align*}
\endgroup
where $\mathcal D:\OO_M\to\cE$ is the derivation determined by $\bracket{\mathcal Df}{e}=\rho(e)f$. We denote this structure by $(\cE,\bracket{-}{-},\rho,\courant{-}{-})$.
\end{dfn}

We use the following algebraic form of Roytenberg's correspondence.

\begin{prop}[Roytenberg correspondence]\label{prop:RoytenbergSevera}
Let $M$ be a smooth complex algebraic variety. Up to isomorphism, Courant algebroids over $M$ correspond bijectively to 2-shifted symplectic varieties $(\mathcal M,\omega,Q)$ of degree $2$ with body $M$. Under this correspondence, $\mathcal O_{\mathcal M}^0=\OO_M$ and $\mathcal O_{\mathcal M}^1\simeq\cE$ as $\OO_M$-modules. Writing $Q=\{\Theta,-\}_{\mathcal M}$ with $\Theta\in\mathcal O(\mathcal M)^3$, the Courant operations are
\begin{align*}
\bracket{e_1}{e_2}=\{e_1,e_2\}_{\mathcal M},\ \rho(e)f=\{\{e,\Theta\}_{\mathcal M},f\}_{\mathcal M},\ 
\courant{e_1}{e_2}=\{\{e_1,\Theta\}_{\mathcal M},e_2\}_{\mathcal M}
\end{align*}
for local sections $e,e_1,e_2\in\cE(U)$ and $f\in\OO_M(U)$.
\end{prop}

\begin{exa}\label{exa:CAPTS}
For $M=\operatorname{Spec}\CC$, the maps $\rho$ and $\mathcal D$ vanish, so a Courant algebroid is precisely a finite-dimensional quadratic Lie algebra $E$, namely a Lie algebra with a nondegenerate invariant symmetric bilinear form. If $E=\fg\oplus\fg^*$ is a decomposition into complementary Lagrangian Lie subalgebras, then $(E,\fg,\fg^*)$ is a Manin triple and $E$ is the Drinfeld double of the corresponding Lie bialgebra. Its graded symplectic model is $E[1]\simeq T^*[2]\fg[1]$, with Hamiltonian given in Example~\ref{exa:LBA}.
\end{exa}

Combining~\cite[Theorem 5.3]{Yam07} with the description of Courant algebroids in terms of 2-shifted symplectic varieties, we obtain the following result.
\begin{prop}[Current PVA associated with a Courant algebroid]\label{prop:PVACA}
Let $(\cE,\bracket{-}{-},\rho,\courant{-}{-})$ be a Courant algebroid over $M$. For each open $U\subseteq M$, take the commutative differential algebra freely generated by $\OO_M(U)$ and $\cE(U)$ with their algebra and module relations, modulo the differential ideal generated by $\partial f-\mathcal Df$ for $f\in\OO_M(U)$. The sheafification $\cV_{\rE}$ of this presheaf carries a unique PVA structure determined by
\begin{align*}
\{f_\lambda g\}_{\rE}=0,\ \{u_\lambda f\}_{\rE}=\rho(u)f,\ \{f_\lambda u\}_{\rE}=-\rho(u)f,\ \{u_\lambda v\}_{\rE}=\courant{u}{v}+\lambda\bracket{u}{v}
\end{align*}
for $f,g\in\OO_M(U)$ and $u,v\in\cE(U)$, extended by sesquilinearity and the Leibniz rules. In the shifted symplectic description, the last formula becomes $\{u_\lambda v\}_{\rE}=\{\{u,\Theta\}_{\mathcal M},v\}_{\mathcal M}+\lambda\{u,v\}_{\mathcal M}.$
\end{prop}

\begin{exa}[Lie bialgebras and their Drinfeld doubles]\label{exa:LBA}
Let $\fg$ be a finite-dimensional complex vector space with basis $\{e_\alpha\}$ and dual basis $\{e^\alpha\}$. On $\mathcal M=T^*[2]\fg[1]$, take degree-$1$ coordinates $q^\alpha,p_\alpha$ corresponding to $e^\alpha,e_\alpha$, respectively, so that $\mathcal O(\mathcal M)^1\simeq\fg^*\oplus\fg$. The canonical symplectic form and Poisson brackets are $\omega=\delta q^\alpha\wedge\delta p_\alpha,\ \{q^\alpha,p_\beta\}_{\mathcal M}=\{p_\beta,q^\alpha\}_{\mathcal M}=\delta_\beta^\alpha,\ 
\{q^\alpha,q^\beta\}_{\mathcal M}=\{p_\alpha,p_\beta\}_{\mathcal M}=0.$
For constants $f_{\alpha\beta}^\gamma=-f_{\beta\alpha}^\gamma$ and $c_\gamma^{\alpha\beta}=-c_\gamma^{\beta\alpha}$, set $\mu=-\frac12f_{\alpha\beta}^\gamma q^\alpha q^\beta p_\gamma,\ \gamma=-\frac12c_\gamma^{\alpha\beta}q^\gamma p_\alpha p_\beta,\ \Theta=\mu+\gamma.$
Since $\{\mu,\mu\}_{\mathcal M}$, $\{\gamma,\gamma\}_{\mathcal M}$ and $\{\mu,\gamma\}_{\mathcal M}$ have distinct $(q,p)$-bidegrees, the master equation is equivalent to $\{\mu,\mu\}_{\mathcal M}=\{\gamma,\gamma\}_{\mathcal M}=\{\mu,\gamma\}_{\mathcal M}=0.$
Define $[-,-]_{\fg}:\wedge^2\fg\to\fg$ and $\delta_{\fg}:\fg\to\wedge^2\fg$ by $[e_\alpha,e_\beta]_{\fg}=f_{\alpha\beta}^\gamma e_\gamma,\ \delta_{\fg}(e_\gamma)=\frac12c_\gamma^{\alpha\beta}e_\alpha\wedge e_\beta,$
so the dual bracket is $[e^\alpha,e^\beta]_{\fg^*}=c_\gamma^{\alpha\beta}e^\gamma$. The first two equations give the Jacobi identities for these brackets, and the third gives $\delta_{\fg}([x,y]_{\fg})=\operatorname{ad}_x\delta_{\fg}(y)-\operatorname{ad}_y\delta_{\fg}(x),\ x,y\in\fg.$
Thus $\{\Theta,\Theta\}_{\mathcal M}=0$ holds precisely when $(\fg,[-,-]_{\fg},\delta_{\fg})$ is a Lie bialgebra. In this case, $(\mathcal M,\omega,\{\Theta,-\}_{\mathcal M})$ is a shifted symplectic variety of degree $2$, and $[u,v]_{\Theta}:=\{\{u,\Theta\}_{\mathcal M},v\}_{\mathcal M}$ on $\mathcal O(\mathcal M)^1$ gives the Drinfeld double $\mathfrak d=\fg\oplus\fg^*$:
\[
[e_\alpha,e_\beta]_{\mathfrak d}=f_{\alpha\beta}^\gamma e_\gamma,\qquad [e^\alpha,e^\beta]_{\mathfrak d}=c_\gamma^{\alpha\beta}e^\gamma,\qquad [e_\alpha,e^\beta]_{\mathfrak d}=c_\alpha^{\beta\gamma}e_\gamma-f_{\alpha\gamma}^\beta e^\gamma.
\]
With the invariant pairing $\bracket{x+\alpha}{y+\beta}=\alpha(y)+\beta(x)$, the triple $(\mathfrak d,\fg,\fg^*)$ is the associated Manin triple.
\end{exa}
\section{$N$-deformations of Courant algebroids}\label{sec:NDCA}

Let $E\to M$ be a Courant algebroid, and let
$(\cM,\omega,\Theta)$ denote the corresponding
degree-$2$ shifted symplectic structure.
Set $\cA:=\calO(\cM),$
so that $\cA^1\simeq\Gamma(E)$, and the degree $-2$ Poisson bracket satisfies $\{u,v\}=\langle u,v\rangle,\  u,v\in\Gamma(E).$
The Courant structure is encoded by a degree-$3$ Hamiltonian $\Theta\in\cA^3$ satisfying $\{\Theta,\Theta\}=0,$
with Dorfman bracket $[u,v]_\Theta:=\{\{u,\Theta\},v\}.$

Let $\cN:E\to E$ be skew-symmetric with respect to $\langle-,-\rangle$. It corresponds to a degree-$2$ function $\wN\in\cA^2$ characterized by $\{u,\wN\}=\cN(u),\ u\in\Gamma(E).$
Set $\Theta_\cN:=\{\wN,\Theta\}.$
Following~\cite{KS11}, we call $\cN$ a \emph{weak deforming tensor for $\Theta$} if $\{\Theta,\{\Theta_\cN,\wN\}\}=0.$
In this case $\Theta_\cN$ is again a Courant Hamiltonian and is compatible with $\Theta$.

From now on, assume that $M$ is a point. Then $E=\fg$ is a finite-dimensional quadratic Lie algebra and $\cM=\fg[1]$. For a skew-symmetric endomorphism $\cN:\fg\to\fg$, define $[x,y]_\cN:=[\cN x,y]+[x,\cN y]-\cN[x,y],$
and its Nijenhuis torsion by $(T_\Theta\cN)(x,y):=[\cN x,\cN y]-\cN[x,y]_\cN.$

\begin{prop}\label{prop:NCYBE}
Let $(\fg,(\cdot,\cdot))$ be a finite-dimensional quadratic Lie algebra, let $\Theta\in\calO(\fg[1])^3$ be the Hamiltonian encoding its Lie bracket, and let $\cN:\fg\to\fg$ be skew-symmetric with corresponding degree-$2$ function $\wN$. Suppose $\cN^2=\kappa\,\Id_\fg,\ \kappa\in\CC.$
Then the following are equivalent:
\begin{enumerate}
\item $\cN$ is Nijenhuis, i.e. $T_\Theta\cN=0$;
\item $\{\{\wN,\Theta\},\wN\}+\kappa\Theta=0;$
\item $[\cN x,\cN y]-\cN\bigl([\cN x,y]+[x,\cN y]\bigr)=-\kappa[x,y],\ x,y\in\fg.$
\end{enumerate}
Hence, for $\kappa=0$, $\cN$ satisfies the operator classical Yang--Baxter equation, while for $\kappa\neq0$ it satisfies the modified classical Yang--Baxter equation with the above convention.
\end{prop}

\begin{proof}
By~\cite[Theorem~3.6]{KS11}, the Nijenhuis torsion is represented by  $\widetilde{T_\Theta\cN}
=-\frac12\left(\{\{\wN,\Theta\},\wN\}+\kappa\Theta\right),$
which proves the equivalence of (1) and (2). Since $\cN^2=\kappa\Id_\fg$,
\[
(T_\Theta\cN)(x,y)
=[\cN x,\cN y]-\cN\bigl([\cN x,y]+[x,\cN y]\bigr)+\kappa[x,y],
\]
which proves the equivalence of (1) and (3).
\end{proof}

\begin{prop}\label{prop:weak-PVA}
Assume that the base variety is a point, and let $\fg$ be the underlying vector space of the Courant algebroid. If $\cN$ is a weak deforming tensor for $\Theta$, then $\Theta_\cN$ defines a Courant structure compatible with $\Theta$. Moreover, the corresponding affine PVA brackets on $\cV:=S(\CC[\partial]\otimes\fg),$
defined on generators by $\{u_\lambda v\}_\Phi:=\{\{u,\Phi\},v\}+\lambda\{u,v\},\
u,v\in\fg,\ \Phi\in\{\Theta,\Theta_\cN\},$
and extended by sesquilinearity and the Leibniz rules, are compatible.
\end{prop}

\begin{proof}
The weak deforming condition and the graded Jacobi identity give $\{\Theta_\cN,\Theta_\cN\}=\{\Theta,\Theta_\cN\}=0$. Fix $\alpha,\beta\in\CC$ and set $\Psi:=\alpha\Theta+\beta\Theta_\cN$ and $k:=\alpha+\beta$. Then $\{\Psi,\Psi\}=0$. Since the base is a point, $[u,v]_\Psi:=\{\{u,\Psi\},v\}$ defines a Lie bracket on $\fg$, and the symmetric bilinear form $\langle u,v\rangle:=\{u,v\}$ is invariant under $[-,-]_\Psi$. For $u,v\in\fg$, write
$
\{u_\lambda v\}:=\alpha\{u_\lambda v\}_\Theta+\beta\{u_\lambda v\}_{\Theta_\cN}=[u,v]_\Psi+k\lambda\langle u,v\rangle.
$
Skew-symmetry follows from the skew-symmetry of $[-,-]_\Psi$ and the symmetry of $\langle-,-\rangle$. For all $u,v,w\in\fg$, we have
\begin{align*}
&\{u_\lambda\{v_\mu w\}\}-\{v_\mu\{u_\lambda w\}\}-\{\{u_\lambda v\}_{\lambda+\mu}w\}\\
&=[u,[v,w]_\Psi]_\Psi-[v,[u,w]_\Psi]_\Psi-[[u,v]_\Psi,w]_\Psi\\
&\quad+k\lambda\bigl(\langle u,[v,w]_\Psi\rangle-\langle[u,v]_\Psi,w\rangle\bigr)
-k\mu\bigl(\langle v,[u,w]_\Psi\rangle+\langle[u,v]_\Psi,w\rangle\bigr)=0
\end{align*}
by the Jacobi identity for $[-,-]_\Psi$ and invariance of the pairing. Since $\cV$ is freely generated by $\fg$ as a commutative differential algebra, sesquilinearity and the Leibniz rules extend these identities to all of $\cV$. Thus every linear combination of the two PVA brackets is a PVA bracket, proving compatibility.
\end{proof}

There are two useful classes of weak deforming tensors arising from the double of a trivial Lie bialgebroid. Let $(A,\mu)$ be a Lie algebroid over $M$, and equip the dual bundle $A^*$ with the trivial Lie algebroid structure. Then $\bigl((A,\mu),(A^*,0)\bigr)$
is a trivial Lie bialgebroid, whose double $E:=A\oplus A^*$ carries the corresponding Courant structure.

Let $\mathcal M=T^*[2]\mathfrak g[1]
\simeq(\mathfrak g\oplus\mathfrak g^*)[1].$
We consider the trivial Lie bialgebra $\bigl((\mathfrak g,\mu),(\mathfrak g^*,0)\bigr),$
so that its doubles
of Lie algebroids is encoded by the degree-$3$ Hamiltonian $\Theta=\mu
=
-\frac12 C_{\alpha\beta}^\gamma
q^\alpha q^\beta p_\gamma .$
An endomorphism $N:\mathfrak g\to\mathfrak g,
\ 
N(e_\beta)=N^\alpha{}_\beta e_\alpha,$
lifts to $\mathcal N=
\begin{pmatrix}
N&0\\
0&-{}^tN
\end{pmatrix}$
on $\mathfrak g\oplus\mathfrak g^*$, and is represented by the
degree-$2$ function $\wN=N^\alpha{}_\beta q^\beta p_\alpha .$
\begin{exa}\label{exa:WDTPN}
Let $\pi\in\Gamma(\wedge^2A)$ and $N:A\to A$ be such that $(\pi,N)$ is a Poisson--Nijenhuis structure on $(A,\mu)$ and $N\circ\pi^\sharp=\pi^\sharp\circ{}^tN.$
Then, by~\cite[Proposition~4.15]{KS11}, the skew-symmetric endomorphism
\[
\cN=
\begin{pmatrix}
N&\pi^\sharp\\
0&-{}^tN
\end{pmatrix}:E\to E
\]
is a weak deforming tensor. Consequently, $\Theta_{\cN}:=\{\widetilde{\cN},\Theta\}$
defines a Courant structure compatible with $\Theta$. In this case, $\Theta_{\cN}
=
\{N,\mu\}+\{\pi,\mu\}
=
\mu_N+\gamma_\pi,$
and it is the doubles
of Lie algebroids of the Lie bialgebroid $\bigl((A,\mu_N),(A^*,\gamma_\pi)\bigr),
\
\gamma_\pi:=\{\pi,\mu\}.$
\end{exa}

\begin{exa}\label{exa:ptcasesl2}
Consider the point case, and let $\fg=\fsl_2$ with basis
$e,h,f$.
Let $e^*,h^*,f^*$ be the dual basis, and equip
$E=\fg\oplus\fg^*$ with the canonical pairing $\langle x+\alpha,y+\beta\rangle=\alpha(y)+\beta(x).$
Take $N=\mathrm{id}_{\fg}$ and $\pi=h\wedge e\in\wedge^2\fg$,
using the contraction convention $\pi^\sharp:\fg^*\longrightarrow\fg,\ 
\pi^\sharp(\alpha)=\alpha(e)h-\alpha(h)e.$
Consider the endomorphism
$
\wN=
\begin{bmatrix}
N & \pi^\sharp\\
0 & -N^t
\end{bmatrix}
\in\End(E).
$
Proposition~\ref{prop:PVACA} associates to the $\wN$-deformed
Courant structure a PVA on
$\mathcal V=\operatorname{Sym}(\CC[\partial]\otimes E)$,
whose bracket is compatible with the undeformed PVA bracket. The brackets among $e,h,f$ are $\{h_\lambda e\}=2e,\ 
\{h_\lambda f\}=-2f,\ 
\{e_\lambda f\}=h.$
Its mixed generator brackets are
\begin{alignat*}{3}
\{e_\lambda e^*\}&=2h^*+\lambda,\qquad
& \{e_\lambda h^*\}&=-f^*,\qquad
& \{e_\lambda f^*\}&=0,\\
\{h_\lambda e^*\}&=2h-2e^*,\qquad
& \{h_\lambda h^*\}&=-2e+\lambda,\qquad
& \{h_\lambda f^*\}&=2f^*,\\
\{f_\lambda e^*\}&=2f,\qquad
& \{f_\lambda h^*\}&=e^*,\qquad
& \{f_\lambda f^*\}&=-2e-2h^*+\lambda,\\
\{h^*_\lambda e^*\}&=-2h^*,\qquad
& \{f^*_\lambda e^*\}&=-2f^*,\qquad
& \{h^*_\lambda f^*\}&=0.
\end{alignat*}
All remaining generator brackets follow by skew-symmetry
or vanish. We denote this Poisson vertex algebra by $\cV_{\wN}$ and the corresponding affine Lie conformal algebra by $\cL_{\wN}$.
\end{exa}

We will compare Examples~\ref{exa:ptcasesl2} and~\ref{exa:RBJDr} in Section~\ref{sec:CRM}.
\section{Chiral r-matrices}\label{sec:CRM}
We now introduce the chiral analogue of a classical \(r\)-matrix. In the terminology following Liberati, it is an \(L\)-invariant solution of the weak conformal classical Yang--Baxter equation. Its geometric interpretation will be developed in forthcoming work of the first author~\cite{Fangpre}.

\begin{dfn}\label{dfn:CCRM}
Let \(L\) be a Lie conformal algebra and let
\(r=\sum_i a_i\otimes b_i\in L\otimes L\). Define
\begin{align*}
\CYB(r):=\sum_{i,j}\Bigg(
\left.[a_i{}_\mu a_j]\otimes b_i\otimes b_j
\right|_{\mu=1\otimes\partial\otimes1}-\left.a_i\otimes[a_j{}_\mu b_i]\otimes b_j
\right|_{\mu=1\otimes1\otimes\partial}-\left.a_i\otimes a_j\otimes[b_j{}_\mu b_i]
\right|_{\mu=1\otimes\partial\otimes1}
\Bigg).
\end{align*}
For notational convenience, we also write
$\llbracket r,r\rrbracket:=\CYB(r)$.
Then \(r\) is called a \emph{chiral classical \(r\)-matrix} if, for every \(x\in L\),
\[
\left.x_\lambda(r+r^{21})\right|_{\lambda=-D_2}=0,
\qquad
\left.x_\lambda\CYB(r)\right|_{\lambda=-D_3}=0,
\]
where \(x_\lambda\) denotes the diagonal conformal adjoint action. We call the second condition the \emph{generalized chiral classical Yang--Baxter equation}; it is also called the weak conformal classical Yang--Baxter equation in the later literature. We say that $r$ satisfies the conformal classical
Yang--Baxter equation (CCYBE) if $\operatorname{CYB}(r)=0
\quad\text{in }L^{\otimes 3}/D_3(L^{\otimes 3}),$
where $D_3=\partial\otimes\Id\otimes\Id
+\Id\otimes\partial\otimes\Id
+\Id\otimes\Id\otimes\partial.$
An element $r\in L\otimes L$ satisfying $r^{21}=-r$
and the CCYBE is called a chiral triangular $r$-matrix.
\end{dfn}
Our definition agrees with that of an $L$-invariant conformal $r$-matrix. We use the term chiral classical $r$-matrix to emphasize the analogy with classical $r$-matrices in ordinary Poisson--Lie theory.
\begin{dfn}\label{dfn:CCLBA}
Let \(L\) be a Lie conformal algebra and let \(r\in L\otimes L\). The triple \(\bigl(L,[\,\cdot_\lambda\cdot\,],r\bigr)\) is called a \emph{coboundary Lie conformal bialgebra} if \(\bigl(L,[\,\cdot_\lambda\cdot\,],dr\bigr)\) is a Lie conformal bialgebra. In this case, \(r\) is called a \emph{coboundary structure} on \(L\).
\end{dfn}

In this section, we compare skew-symmetric classical r-matrices with classical r-matrices in the chiral setting.
\begin{thm}\label{thm:CCYBTOOCYB}
Let \(L=\Cur(\fg)\), where \(\fg\) is a quadratic Lie algebra with a nondegenerate symmetric invariant bilinear form \((\cdot,\cdot)\). Denote by \((\cdot,\cdot)_\lambda\) the induced conformal bilinear form on \(L\), as in Example~\ref{exa:CBL}.
Let $r=\sum_{i,n}C_{i,n}a_i\otimes b_i^{(n)}\in L\otimes L,\ b_i^{(n)}:=\partial^n b_i,$
be a skew-symmetric solution of the conformal classical Yang--Baxter equation. For a formal parameter \(\rho\), define the \(\CC\)-linear map \(R_\rho:L\to L[\rho]\) by $R_\rho(x):=\sum_{i,n}C_{i,n}(b_i^{(n)},x)_\rho a_i,\ x\in L.$
Equivalently, by conformal symmetry, $R_\rho(x)=\sum_{i,n}C_{i,n}(x,b_i^{(n)})_{-\rho}a_i.$
Set $R(x):=R_\partial(x)=\sum_{i,n}C_{i,n}(b_i^{(n)},x)_\partial a_i,$
where the substituted \(\partial\) acts on \(a_i\). Then \(R\) is a \(\CC[\partial]\)-module endomorphism of \(L\) satisfying the operator classical Yang--Baxter equation for Lie conformal algebras~\cite[Eq.~8]{Fang25}, namely, $[R(x)_\lambda R(y)]
=R([R(x)_\lambda y]+[x_\lambda R(y)])$ for all $x,y$ in $L.$
\end{thm}
\begin{proof}
Write \(b_i^{(n)}:=\partial^n b_i\), and for \(P(\rho)=\sum_{k\geq0}p_k\rho^k\in\mathbb C[\rho]\) set \(P(\partial)u:=\sum_{k\geq0}p_k\partial^ku\).

Let $r=\sum_{i,n}C_{i,n}a_i\otimes b_i^{(n)}\in L\otimes L$
be skew-symmetric, \(q,x,y,z\in L\), and define
\[
R_\rho(x):=\sum_{i,n}C_{i,n}( b_i^{(n)},x)_\rho a_i\in L[\rho],\qquad R(x):=R_\partial(x)=\sum_{i,n}C_{i,n}(b_i^{(n)},x)_\partial a_i,
\]
where the substituted \(\partial\) acts on \(a_i\). By sesquilinearity,  $R(\partial^kx)=\sum_{i,n}C_{i,n}\partial^k\bigl(( b_i^{(n)},x)_\partial a_i\bigr)=\partial^kR(x),$
so \(R\in\operatorname{End}_{\mathbb C[\partial]}(L)\).

The skew-symmetry \(r^{21}=-r\) gives
\[
\sum_{i,n}C_{i,n}(a_i,q)_\rho b_i^{(n)}=-R_\rho(q),
\tag{1}
\]
and also
\[
(R(q),x)_{-\lambda}=\sum_{i,n}C_{i,n}( b_i^{(n)},q)_\lambda( a_i,x)_{-\lambda},
\tag{2}
\]
Here we use sesquilinearity and apply the pairing
$(\,\cdot\,,q)_\lambda$ to the identity $\sum_{i,n}C_{i,n}(a_i,x)_{-\lambda}b_i^{(n)}
=-R_{-\lambda}(x)
=-\sum_{i,n}C_{i,n}(b_i^{(n)},x)_{-\lambda}a_i.$

In addition, since $r^{21}=-r$ we have
\[
( R_{-\rho}(q),z)_\rho=-( q,R(z))_\rho,\qquad ( R_\rho(q),z)_{-\rho}=-( R_{-\rho}(z),q)_\rho.
\tag{3}
\]

Write $\llbracket r,r\rrbracket=\sum_{i,j,n,m}C_{i,n}C_{j,m}\left(T^{(1)}_{i,j,n,m}-T^{(2)}_{i,j,n,m}-T^{(3)}_{i,j,n,m}\right),$
where
\[
\begin{gathered}
T^{(1)}_{i,j,n,m}
:=\left.[a_i{}_\mu a_j]\otimes b_i^{(n)}\otimes b_j^{(m)}
\right|_{\mu=\partial_2},
\qquad
T^{(2)}_{i,j,n,m}
:=\left.a_i\otimes[a_j{}_\mu b_i^{(n)}]\otimes b_j^{(m)}
\right|_{\mu=\partial_3},
\\
T^{(3)}_{i,j,n,m}
:=\left.a_i\otimes a_j\otimes[b_j^{(m)}{}_\mu b_i^{(n)}]
\right|_{\mu=\partial_2}.
\end{gathered}
\]
where \(\partial_1:=\partial\otimes1\otimes1\), \(\partial_2:=1\otimes\partial\otimes1\), and \(\partial_3:=1\otimes1\otimes\partial\), with \(1\) denoting the identity map on \(L\).
For \(A,B,C,x,y,z\in L\), set $( A\otimes B\otimes C,x\otimes y\otimes z)_{\alpha,\beta,\gamma}:=( A,x)_\alpha(B,y)_\beta( C,z)_\gamma.$
Take $\alpha=-\lambda,\  \beta=\lambda-\eta,\  \gamma=\eta,\ \nu:=\eta-\lambda,$
so that \(\beta=-\nu\) and \(\alpha+\beta+\gamma=0\).

Let \(E_1,E_2,E_3\) denote the contributions of \(T^{(1)},-T^{(2)},-T^{(3)}\), respectively. Since \(\partial_2\) becomes \(-\beta=\nu\),
\[
\begin{aligned}
E_1&=\sum_{i,j,n,m}C_{i,n}C_{j,m}([a_i{}_\nu a_j],x)_{-\lambda}( b_i^{(n)},y)_{-\nu}( b_j^{(m)},z)_\eta\\
&=([R_{-\nu}(y){}_\nu R_\eta(z)],x)_{-\lambda}=([R(y){}_\nu R(z)],x)_{-\lambda}.
\end{aligned}
\tag{4}
\]
The last equality holds after applying the outer pairing: it follows from $[P(\partial)a{}_\nu Q(\partial)b]=P(-\nu)Q(\nu+\partial)[a{}_\nu b]$
and the fact that the remaining \(\partial\) becomes \(\lambda\), so that \(\nu+\lambda=\eta\).

For the second term, \(\partial_3\) becomes \(-\eta\). Put \(w_{i,n}:=[b_i^{(n)}{}_{-\nu-\partial}y]\). By invariance and (1)--(3),
\[
\begin{aligned}
E_2&=-\sum_{i,j,n,m}C_{i,n}C_{j,m}( a_i,x)_{-\lambda}([a_j{}_{-\eta}b_i^{(n)}],y)_{-\nu}( b_j^{(m)},z)_\eta\\
&=\sum_{i,n}C_{i,n}(a_i,x)_{-\lambda}( R_{-\eta}(w_{i,n}),z)_\eta\\
&=-\sum_{i,n}C_{i,n}(a_i,x)_{-\lambda}(w_{i,n},R(z))_\eta.
\end{aligned}
\]
Inside the indicated pairings, sesquilinearity and invariance give
\[
( w_{i,n},R(z))_\eta=([b_i^{(n)}{}_\lambda y],R(z))_\eta=( b_i^{(n)},[y{}_\nu R(z)])_\lambda.
\]
Therefore, by (2),
\[
E_2=-( R([y{}_\nu R(z)]),x)_{-\lambda}.
\tag{5}
\]

For the third term, \(\partial_2\) again becomes \(\nu\). Put \(W_{i,n}:=[b_i^{(n)}{}_{\eta-\partial}z]\). By invariance and (3),
\[
\begin{aligned}
E_3&=-\sum_{i,j,n,m}C_{i,n}C_{j,m}( a_i,x)_{-\lambda}( a_j,y)_{-\nu}([b_j^{(m)}{}_\nu b_i^{(n)}],z)_\eta\\
&=-\sum_{i,n}C_{i,n}(a_i,x)_{-\lambda}( R_\nu(W_{i,n}),y)_{-\nu}\\
&=\sum_{i,n}C_{i,n}(a_i,x)_{-\lambda}( R_{-\nu}(y),W_{i,n})_\nu. 
\end{aligned}\tag{6}
\]
Since the \(\partial\) in \(W_{i,n}\) occurs in the second argument, it becomes \(\nu\), and hence \(\eta-\partial\) becomes \(\eta-\nu=\lambda\). Moreover, for \(A,b,z\in L\), $( A,[b{}_\lambda z])_\nu=([b{}_\lambda z],A)_{-\nu}=( b,[z{}_{-\nu-\partial}A])_\lambda=-( b,[A{}_\nu z])_\lambda.$
Thus, using \([R_{-\nu}(y){}_\nu z]=[R(y){}_\nu z]\), $E_3=-\sum_{i,n}C_{i,n}(a_i,x)_{-\lambda}(b_i^{(n)},[R(y){}_\nu z])_\lambda
=-(R([R(y){}_\nu z]),x)_{-\lambda}.$

Combining (4)--(6), we obtain
\[
\begin{aligned}
&(\llbracket r,r\rrbracket,x\otimes y\otimes z)_{-\lambda,\lambda-\eta,\eta}=\left([R(y){}_\nu R(z)]-R\bigl([R(y){}_\nu z]+[y{}_\nu R(z)]\bigr),x\right)_{-\lambda}.
\end{aligned}
\tag{7}
\]

Let \(D_3:=\partial_1+\partial_2+\partial_3\). Since $\llbracket r,r\rrbracket=0\pmod{D_3(L^{\otimes3})},$
there exists \(Q\in L^{\otimes3}\) such that \(\llbracket r,r\rrbracket=D_3Q\). But $(D_3Q,x\otimes y\otimes z)_{\alpha,\beta,\gamma}=-(\alpha+\beta+\gamma)( Q,x\otimes y\otimes z)_{\alpha,\beta,\gamma}=0.$
Hence the right hand side of (7) vanishes. Setting \(\eta=\lambda+\nu\), define $F_\nu(y,z):=[R(y){}_\nu R(z)]-R\bigl([R(y){}_\nu z]+[y{}_\nu R(z)]\bigr)\in L[\nu].$
Then \(( F_\nu(y,z),x)_{-\lambda}=0\) for all \(x,y,z\in L\). By the nondegenerate of conformal bilinear form, $F_\nu(y,z)=0$ for all $y,z\in L$, hence $F_\nu=0.$
Therefore $[R(y){}_\nu R(z)]=R\bigl([R(y){}_\nu z]+[y{}_\nu R(z)]\bigr),$
so \(R\) satisfies the operator conformal classical Yang--Baxter equation.
\end{proof}
\setcounter{equation}{7}
In the following examples, we use the normalized invariant bilinear form on \(\fsl_2\) given by \(\kappa(e,f)=\kappa(f,e)=1\) and \(\kappa(h,h)=2\).

\begin{exa}\label{exa:RBJDr}
Consider \(L=\Cur(\fsl_2)\), with the convention \(h\wedge e:=h\otimes e-e\otimes h\). For \(c\in\CC\), define
\[
r_c=h\wedge e+\frac{c}{2}\bigl(\partial e\otimes e-e\otimes\partial e\bigr)\in L\otimes L.
\]
This is a skew-symmetric solution of the conformal classical Yang--Baxter equation. By Theorem~\ref{thm:CCYBTOOCYB}, it determines a classical \(R\)-matrix $R_c(x)=\kappa_\partial(e,x)h-\kappa_\partial(h,x)e+c\kappa_\partial(e,x)\partial e.$
By the definition of the \(R\)-twisted bracket in~\cite{Fang25}, the nonzero \(\lambda\)-brackets of \(L_{R_c}\) are
\begin{align*}
[e{}_\lambda f]_{R_c}&=-2e,\qquad
[h{}_\lambda f]_{R_c}=-2h+2c\lambda e+2c\partial e,\\
[f{}_\lambda f]_{R_c}&=-c(2\lambda+\partial)h,\qquad
[f{}_\lambda e]_{R_c}=2e,\\
[f{}_\lambda h]_{R_c}&=2h+2c\lambda e.
\end{align*}
\end{exa}
\begin{rmk}
We compare Examples~\ref{exa:ptcasesl2} and~\ref{exa:RBJDr}. Let $\cL$ denote the affine LCA in Example~\ref{exa:ptcasesl2}. For $c\in\CC^\times$, set $E:=f^*$, $H:=2h^*+c\partial f^*$, and $F:=e^*-c\partial h^*-\frac{c^2}{2}\partial^2f^*$, and define $\cL_-:=\CC[\partial]E\oplus\CC[\partial]H\oplus\CC[\partial]F\subset\cL$. Since this change of generators is invertible over $\CC[\partial]$, we have $\cL_-=\CC[\partial]\fsl_2^*$; in particular, this submodule is independent of $c$. A direct computation of the inherited $\lambda$-brackets shows that $\cL_-$ is a Lie conformal subalgebra of $\cL$ and that the $\CC[\partial]$-linear map
$
\Phi_c:L_{R_c}\longrightarrow\cL_-,\  e\longmapsto E,\  h\longmapsto H,\  f\longmapsto F,
$
is an isomorphism of Lie conformal algebras.
\end{rmk}

\begin{prop}
Set $\cL_+:=\CC[\partial]\fsl_2$. Then
$\cL=\cL_+\oplus\cL_-\oplus\CC K$
as a $\CC[\partial]$-module, where $K$ is central and $\partial K=0$. Equip $\cD:=\cL_+\oplus\cL_-$ with the $\lambda$-bracket obtained from that of $\cL$ by omitting the central terms. Then $\cL$ is a one-dimensional central extension of the Lie conformal algebra $\cD$.

The canonical pairing on $\fsl_2\oplus\fsl_2^*$ extends to a nondegenerate symmetric invariant conformal bilinear form on $\cD$, with respect to which $(\cD,\cL_+,\cL_-)$ is a conformal Manin triple, in the sense of ~\cite[Definition 3.8]{Lib08}.
Thus $\cD$ is the corresponding conformal double.

Moreover, the associated universal affine vertex algebra $U_{\wN}$ at level one is a vertex operator algebra of central charge $6$, with conformal vector
$
\omega=:e(e^*-h):+:h(h^*+e):+:ff^*:-4:(e^*-h)f^*:-4:(h^*+e)(h^*+e):.$
This vertex algebra embeds into $\cD_{\SL_2,0}^{\ch}$, the vertex algebra of global chiral differential operators on $\SL_2$ at level zero.
\end{prop}
\begin{proof}
It is straightforward to check that $\omega$ is a conformal
vector of central charge $6$.
With respect to the deformed bracket and the canonical pairing, $I:=\Span\{e^*-h,h^*+e,f^*\}$ is an abelian Lagrangian ideal of $\fsl_2\oplus\fsl_2^*$, and $\fsl_2\oplus\fsl_2^*=\fsl_2\oplus I$. Hence $\mathfrak d:=\fsl_2\oplus I$ is isomorphic to the cotangent Lie algebra $T^*\fsl_2$ as a quadratic Lie algebra. By \cite[Corollary~4.17]{AN25}, we therefore obtain an embedding $V_1(\mathfrak d)\hookrightarrow\Gamma(G,\cD^\ch_{G,0})$, where $G=\mathrm{SL}_2(\mathbb C)$. This statement extends to any finite-dimensional complex Lie algebra $\fg$. Let $\fd:=\fg\oplus\fg^*$ be the Drinfeld double of $\fg$ with trivial cobracket, and let $\pi\in\wedge^2\fg$ satisfy the classical Yang--Baxter equation. Set $\cN:=\begin{bmatrix}
\Id_\fg & \pi^\#\\
0 & -\Id_{\fg^*}
\end{bmatrix},
\  V:=\operatorname{graph}(-\pi^\#).$
Then $V\cong\fg^*$ is an abelian Lagrangian ideal of $\fd_\cN$, and the decomposition $\fd_\cN=\fg\oplus V$ identifies $\fd_\cN$ with the cotangent Lie algebra $T^*\fg$ as a quadratic Lie algebra. If moreover $\fg=\operatorname{Lie}(G)$ for a connected complex reductive algebraic group $G$, then \cite[Corollary~4.17]{AN25} yields an embedding $V_1(\fd_\cN,B)\hookrightarrow\Gamma(G,\cD^\ch_{G,0}),$
where $B$ is the canonical pairing and $V_1(\fd_\cN,B)$ denotes the corresponding universal affine vertex algebra.
\end{proof}
The first author will study the conformal double structure further in his forthcoming paper.

\begin{exa}\label{exa:RBJCa}
Let \(L=\Cur(\fsl_2)\). A family of classical \(r\)-matrices is given by~\cite{GK23} $r^{(h)}_p=\frac{1}{4}\bigl(p(\partial)h\otimes h-h\otimes p(\partial)h\bigr),$
where \(p(x)\) is a nonzero odd polynomial, i.e., \(p(-x)=-p(x)\). Here the subscript \(h\) indicates the Cartan direction
\(\mathbb C h\), while the dependence on the fixed polynomial
\(p\) is suppressed in the notation. By Theorem~\ref{thm:CCYBTOOCYB}, the corresponding \(R\)-matrix is given by $R_h(e)=R_h(f)=0,\ R_h(h)=p(\partial)h.$
The nonzero \(\lambda\)-brackets of the \(R\)-deformed LCA \(L_{R_h}\) are $[h_\lambda f]_{R_h}=2p(\lambda)f,\ 
[e_\lambda h]_{R_h}=-2p(\lambda+\partial)e,
[f_\lambda h]_{R_h}=2p(\lambda+\partial)f,
[h_\lambda e]_{R_h}=-2p(\lambda)e.$

\end{exa}
\section{Applications to integrable systems}\label{sec:APPIS}

In~\cite{GK23}, the authors classified the solutions of the conformal classical Yang--Baxter equation on $\Cur(\fsl_2)$, while in~\cite{GKRB23}, they classified Rota--Baxter operators on $\Cur(\fsl_2)$ and determined which of them arise from such solutions. Using these results together with Theorem~\ref{thm:CCYBTOOCYB}, we construct several examples of integrable systems via the generalized AKS scheme (GAKS)~\cite{Fang25,Fang26} and the Lenard--Magri scheme.
\begin{exa}\label{exa:HRBJDr}
With \(u=(f,h,e)^T\) and \(H_{ij}(\lambda):=[u_i{}_\lambda u_j]_{R_c}\), the matrix differential operator associated with \(R_c\) in Example~\ref{exa:RBJDr} decomposes as $H_{R_c}=H^{(0)}+cH^{(1)},$
where
\[
H^{(0)}=\begin{pmatrix}
0&2h&2e\\
-2h&0&0\\
-2e&0&0
\end{pmatrix},\qquad
H^{(1)}=\begin{pmatrix}
-2h\partial-h^{(1)}&2e\partial&0\\
2e\partial+2e^{(1)}&0&0\\
0&0&0
\end{pmatrix}.
\]
Since \(H_{R_c}\) is Hamiltonian for every \(c\in\CC\), \(H^{(0)}\) and \(H^{(1)}\) form a compatible Hamiltonian pair. For column variational derivatives, set \(P_i:=(H^{(i)})^{\mathsf T}\), \(i=0,1\).

Fix \(\alpha\in\CC\) and set $d:=\alpha e-h,\ \widetilde{\mathcal V}:=\mathcal V[d^{-1}],\ \mathcal F:=\widetilde{\mathcal V}/\partial\widetilde{\mathcal V}.$
For every \(n\geq0\), define
\[
\mathcal H_n:=\int h_n,\qquad h_n:=\frac{e(e^{(1)})^n}{d^{n+1}}.
\]
Then
\[
P_0\frac{\delta\mathcal H_0}{\delta u}=0,\qquad
P_1\frac{\delta\mathcal H_n}{\delta u}
=P_0\frac{\delta\mathcal H_{n+1}}{\delta u}\quad(n\geq0).
\]
Consequently, the corresponding hierarchy is
\[
\frac{\partial f}{\partial t_n}
=2(n+1)\partial\left(\frac{e^2(e^{(1)})^n}{(\alpha e-h)^{n+2}}\right),\qquad
\frac{\partial h}{\partial t_n}=\frac{\partial e}{\partial t_n}=0.
\]
\end{exa}
\begin{exa}\label{exa:HRBJCa}
Let \(R_h\) be the \(R\)-matrix constructed in Example~\ref{exa:RBJCa}. Following~\cite{Fang25}, let \(\mathcal V:=\Sym(L)\), where \(L=\Cur(\fsl_2)\), be the PVA induced by the \(R_h\)-deformed bracket. With respect to \(u=(f,h,e)^T\) and the convention \(H_{ij}(\lambda)=[u_i{}_\lambda u_j]_{R_h}\), its Hamiltonian operator is
\begin{equation}\label{eq:HRh}
H_{R_h}(\partial)=
\begin{bmatrix}
0&2p(\partial)\circ f&0\\
2f\,p(\partial)&0&-2e\,p(\partial)\\
0&-2p(\partial)\circ e&0
\end{bmatrix}.
\end{equation}
Consider the localization of $\cV$ at $e$ and denote the localized differential algebra by $\cV_e$ Then $\cV_e=\cV[e^{-1}].$
 In particular,
$\{a_\lambda e^{-1}\}_R=-e^{-2}\{a_\lambda e\}_R,\  a\in\mathcal V_e.$
Set \(\rho:=ef\). It is clear, $\rho\in\fz(\cV_e),$ where $\fz(\cV_e)$ is the Poisson center of $\cV_e.$ Therefore, as PVAs,
$
\mathcal V_e\cong
\CC\bigl[h^{(n)},e^{(n)}\mid n\geq0\bigr][e^{-1}]
\otimes
\CC\bigl[\rho^{(n)}\mid n\geq0\bigr],
$
where the second factor has zero \(\lambda\)-bracket.
\end{exa}
\section{Applications to vertex algebras}\label{sec:APVAS}

Following~\cite{Fang25}, we apply the preceding results to vertex algebras. In particular, by combining our $r$-matrix construction with the classifications in~\cite{GK23,GKRB23}, we construct several examples of vertex algebras associated with the resulting classical r-matrices.
Also, this deformed Lie conformal algebra can be quantized by the universal enveloping vertex algebra~\cite{Fang25,Primc99,DSK05}.
\begin{prop}\label{prop:PRC}
Let $\fg$ be a Lie algebra equipped with an invariant symmetric bilinear
form $(\,\cdot\mid\cdot\,)$, and let $L=\Cur(\fg)$. Suppose that
$R\in\End_{\CC[\partial]}(L)$ satisfies classical CYBE of $L$ and denote the $R$-deformed LCA by $L_R.$ Let $\hat L$ be the affine central extension of LCA $L,$ whose $\lambda$-bracket
on $\fg$ is given by $[a_\lambda b]^{\wedge}
=
[a,b]+k\lambda(a\mid b)\vac.$
Denote the corresponding $2$-cocycle by $\omega_{k,\lambda}(a,b):=k\lambda(a\mid b).$
Then the pullback $(R^*\omega_k)_\lambda(a,b)
:=
\omega_{k,\lambda}(R(a),R(b))$
is a $2$-cocycle of $L_R$. Consequently, it defines a central extension
$\widehat L_R=L_R\oplus\CC\vac$
with $\lambda$-bracket $[a_\lambda b]^{\wedge}_R
=
[a_\lambda b]_R
+
(R^*\omega_k)_\lambda(a,b)\vac.$
\end{prop}

\begin{proof}
The classical Yang--Baxter equation implies that $R([X_\lambda Y]_R)=[R(X)_\lambda R(Y)],$
so that $R:L_R\to L$ is a homomorphism of Lie conformal algebras.
We need to verify that $(R^*\omega_k)_\lambda$ is a $2$-cocycle, i.e., $(R^*\omega_k)_\lambda(X,[Y_\mu Z]_R)
-
(R^*\omega_k)_\mu(Y,[X_\lambda Z]_R)
=
(R^*\omega_k)_{\lambda+\mu}([X_\lambda Y]_R,Z)$
for all $X,Y,Z\in L_R$. This follows immediately from the fact that
$R$ is a homomorphism of Lie conformal algebras and that $\omega_k$
is a $2$-cocycle of $L$.
\end{proof}
\begin{thm}\label{thm:EEXTsl2}
Let $\fg=\fsl_2$, let $L=\Cur(\fg)$, and consider $\hat L_R$,
where $R$ is the classical $R$-matrix in Example~\ref{exa:RBJCa}.
There exists a 2-cocycle $\eta$, linearly independent of $R^*\omega_k,\ for\  k\neq0$, given by $\eta_\lambda(f,e)=1,\ \eta_\lambda(e,f)=-1,$
and zero on all other pairs of generators.
We denote by $\widetilde L_R$ the central extension associated with
the cocycle $R^*\omega_k+\kappa\eta$.
Its nontrivial $\lambda$-brackets are $[h_\lambda e]_R^\sim
=
[h_\lambda e]_R,\ 
[h_\lambda f]_R^\sim
=
[h_\lambda f]_R,
\ 
[h_\lambda h]_R^\sim
=
-2k\lambda p(\lambda)^2K,\ 
[f_\lambda e]_R^\sim
=\kappa K.$
Assume that $\kappa\neq0$, and set $V_{p,k,\kappa}
:=
U(\widetilde L_R)/(K-\mathbf1).$
Let $\mathcal S_\kappa$ be the rank one Weyl vertex algebra generated by
$\beta,\gamma$ $[\gamma_\lambda\beta]
=
\kappa\mathbf1.$
Let $\mathcal H_{2(k+2)}$ be the rank one Heisenberg vertex algebra
generated by $\varphi$ with $[\varphi_\lambda\varphi]
=
2(k+2)\lambda\mathbf1,$
and let $\mathcal H^{(p)}_{2(k+2)}
:=
\bigl\langle p(\partial)\varphi\bigr\rangle_{\mathrm{VA}}
\subset
\mathcal H_{2(k+2)}$
be the vertex subalgebra generated by $p(\partial)\varphi$.
For $k=-2$, we take $\mathcal H_0=\mathbb C[\partial^n\varphi\mid n\geq0]$ with its standard commutative vertex algebra structure.
Then $V_{p,k,\kappa}
\cong
\mathcal S_\kappa
\otimes
\mathcal H^{(p)}_{2(k+2)}.$ In particular, there is an injective free field realization $V_{p,k,\kappa}
\hookrightarrow
\mathcal S_\kappa
\otimes
\mathcal H_{2(k+2)}.$
\end{thm}
\begin{proof}
A direct calculation shows that $\eta$ is a $2$-cocycle and, for $k\ne0$, is linearly independent of $R^*\omega_k$. Hence $R^*\omega_k+\kappa\eta$ defines a central extension for every $k\in\CC$.

Set $u:=:\!\beta\gamma\!:$ and $t:=p(\partial)\varphi$. The Wick formula, sesquilinearity, and the oddness of $p$ give
\[
[u_\lambda\beta]=\kappa\beta,\qquad [u_\lambda\gamma]=-\kappa\gamma,\qquad [u_\lambda u]=-\kappa^2\lambda,\qquad [t_\lambda t]=-2(k+2)\lambda p(\lambda)^2.
\]
Define $\Phi:V_{p,k,\kappa}\to\mathcal S_\kappa\otimes\mathcal H^{(p)}_{2(k+2)}$ by $\Phi(e)=\beta,\ \Phi(f)=\gamma,\ \Phi(h)=t+\frac{2}{\kappa}p(\partial)u.$
Since the Weyl and Heisenberg factors commute, the preceding identities imply
\[
[\Phi(h)_\lambda\Phi(e)]=-2p(\lambda)\beta,\qquad [\Phi(h)_\lambda\Phi(f)]=2p(\lambda)\gamma,\qquad [\Phi(h)_\lambda\Phi(h)]=-2k\lambda p(\lambda)^2,\qquad [\Phi(f)_\lambda\Phi(e)]=\kappa.
\]
Thus $\Phi$ preserves the defining brackets. Set $a:=h-\frac{2}{\kappa}p(\partial):\!ef\!:.$
Then $\Phi(a)=t$. A direct application of the Wick formula gives $[f_\lambda e]=\kappa,\
[a_\lambda e]=[a_\lambda f]=0,\
[a_\lambda a]=-2(k+2)\lambda p(\lambda)^2.$
Since $\mathcal H_{2(k+2)}$ is freely generated by $\varphi$ and
$p\neq0$, the $\mathbb C[\partial]$-submodule
$\mathbb C[\partial]p(\partial)\varphi$ is free of rank one. By the PBW
basis of the Heisenberg vertex algebra, the normally ordered monomials in
the derivatives of $p(\partial)\varphi$ are linearly independent. Hence
$\mathcal H^{(p)}_{2(k+2)}$ is freely generated by
$t=p(\partial)\varphi$.
Therefore the PBW presentations give a vertex algebra homomorphism $\Psi:\mathcal S_\kappa\otimes
\mathcal H^{(p)}_{2(k+2)}
\longrightarrow V_{p,k,\kappa}$
defined by $\Psi(\beta)=e,\ 
\Psi(\gamma)=f,\ 
\Psi(t)=a.$
The maps $\Phi$ and $\Psi$ are inverse on generators, hence they are
inverse vertex algebra homomorphisms. Therefore $\Phi$ is an isomorphism.
\end{proof}
We now relate the vertex algebra in Theorem~\ref{thm:EEXTsl2}
to the work of~\cite{CL22,ACKL17}. Following these works, we use the notation $\OO_{\ev}(n,d),\  \cS_{\ev}(n,d),\ 
\OO_{\odd}(n,d),\  \cS_{\odd}(n,d)$
for the standard free field algebras. Here $\OO$ and $\cS$ denote the
orthogonal and symplectic types of the bilinear form determined by the
leading OPE coefficients, while $\ev$ and $\odd$ indicate the parity of
the generating fields.

For even $d$, the algebra $\OO_{\ev}(n,d)$ is generated by $n$ even fields
$a^1,\ldots,a^n$ of conformal weight $d/2$, with $a^i(z)a^j(w)\sim\frac{\delta_{ij}}{(z-w)^d}.$
For odd $d$, the algebra $\cS_{\ev}(n,d)$ is generated by $2n$ even fields
$a^1,\ldots,a^n,b^1,\ldots,b^n$ of conformal weight $d/2$, whose only
nonzero singular OPEs are $a^i(z)b^j(w)\sim\frac{\delta_{ij}}{(z-w)^d},
\ 
b^i(z)a^j(w)\sim-\frac{\delta_{ij}}{(z-w)^d}.$
\begin{prop}\label{prop:EEXTsl2}
Let $p(x)=x$ and set $A:=h-\frac{2}{\kap}\partial(:ef:)$. If $k\neq-2$ and $\kap\neq0$, then
\[
V_{p,k,\kap}\cong\cS_{\ev}(1,1)\otimes\langle A\rangle_{\mathrm{VA}}\cong\cS_{\ev}(1,1)\otimes\OO_{\ev}(1,4).
\]
Moreover, the vertex subalgebra $V'_{p,k,\kap}:=\langle\partial e,\partial f\rangle_{\mathrm{VA}}\subset V_{p,k,\kap}$ satisfies $V'_{p,k,\kap}\cong\cS_{\ev}(1,3)$, and
\[
\cW_{\free}(\fsl_3,f_{\min})\cong\OO_{\ev}(1,2)\otimes V'_{p,k,\kap}\otimes\langle A\rangle_{\mathrm{VA}}\cong\OO_{\ev}(1,2)\otimes\cS_{\ev}(1,3)\otimes\OO_{\ev}(1,4).
\]
\end{prop}

\begin{proof}
Write $\cS_{\ev}(1,1)=\langle\beta,\gamma\rangle_{\mathrm{VA}}$ and $\OO_{\ev}(1,4)=\langle a\rangle_{\mathrm{VA}}$, where $[\gamma_\lambda\beta]=1$ and $[a_\lambda a]=\frac{1}{6}\lambda^3$. By the preceding theorem, $[f_\lambda e]=\kap,\  [A_\lambda e]=[A_\lambda f]=0,\ [A_\lambda A]=-2(k+2)\lambda^3.$
Choose $c\in\mathbb C^\times$ with $c^2=-12(k+2)$. Since $h=A+\frac{2}{\kap}\partial(:ef:)$, the assignments $e\mapsto\beta,\  f\mapsto\kap\gamma,\ A\mapsto ca$
define an isomorphism $V_{p,k,\kap}\cong\cS_{\ev}(1,1)\otimes\OO_{\ev}(1,4)$, whose inverse sends $\beta\mapsto e$, $\gamma\mapsto\kap^{-1}f$, and $a\mapsto c^{-1}A$. In particular, $\langle A\rangle_{\mathrm{VA}}\cong\OO_{\ev}(1,4)$. Moreover, sesquilinearity gives $[\partial e_\lambda\partial f]=\kap\lambda^2$. Since the standard generators $b,c$ of $\cS_{\ev}(1,3)$ satisfy $[b_\lambda c]=\frac{1}{2}\lambda^2$, after rescaling one generator the PBW theorem gives $V'_{p,k,\kap}\cong\cS_{\ev}(1,3)$. The last assertion follows from $\cW_{\free}(\fsl_3,f_{\min})\cong\OO_{\ev}(1,2)\otimes\cS_{\ev}(1,3)\otimes\OO_{\ev}(1,4).$
\end{proof}

\section{The Virasoro obstruction via shifted symplectic chiralization}\label{sec:OQVOS}

The preceding results, together with~\cite{Fang25,Fang26,FangPhD}, motivate our study of 2-shifted symplectic chiralization and the quantization of Lie conformal bialgebra structures. Here, we verify that our 2-shifted symplectic chiralization detects precisely the same obstruction to the Lie conformal bialgebra compatibility condition in the Virasoro case, up to the sign convention for the variational bracket.

\begin{prop}\label{prop:virob}
Let $\Vir=\mathbb{C}[\partial]L$ be the Virasoro conformal algebra, and let $\Vir^*_{\ch}=\mathbb{C}[\partial]L^*$
be its conformal dual, normalized by $(L^*)_\lambda(L)=1,
\ 
(L^*)_\lambda(\partial^n L)=\lambda^n (L^*)_\lambda(L).$
Define a cobracket on $\Vir$ by $\delta(L)=L^{(1)}\otimes L-L\otimes L^{(1)}.$
Then the induced Lie conformal algebra on the conformal dual satisfies $(\Vir^*_{\ch},\delta^*_{\ch})\cong (\Vir,[\,{}_\lambda{}\,]).$
However, this cobracket does not make $\Vir$ into a Lie conformal
bialgebra. More precisely, the Virasoro $\lambda$-bracket and the
cobracket $\delta$ are not compatible; instead, one has
\[
L_\lambda\delta(L)-L_{-\lambda-D_2}\delta(L)
=
5\,\delta([L_\lambda L]).\]
\end{prop}
\begin{proof}
This follows directly from \cite[Proposition~2.13(a) and Example~2.15]{Lib08}.
\end{proof}

Recall that in Poisson--Lie theory, if $\mu$ denotes the Lie bracket on
$\fg$ and $\gamma$ denotes the cobracket on $\fg$, then the compatibility
between $\mu$ and $\gamma$ is encoded by the condition $\{\mu,\gamma\}=0 .$
Motivated by~\cite{Fang26,FangPhD} and by the idea of
$n$-shifted symplectic chiralization proposed by the first author, we consider $X=\mathbb{A}^1,\
Y:=T^*[1]\mathbb{A}^1,\
Z:=\mathbb{A}^1[1]\oplus(\mathbb{A}^1)^*[1],$
together with their arc spaces. In this setting, the Virasoro
$\lambda$-bracket and the cobracket on its conformal dual are represented
by the following two degree three local functionals: $\mu=\int \xi\theta\theta^{(1)},\ and\ 
\gamma=\int \theta\xi\xi^{(1)}.$
Here $\xi$ is the degree one coordinate on $\mathbb{A}^1[1]$, and
$\theta$ is the degree one fiber coordinate on $(\mathbb{A}^1)^*[1]$.
Equivalently, $\mu$ and $\gamma$ are local functionals on $Z_\infty$,
i.e. elements of $\mathcal{F}_{Z_\infty}:=\mathcal{O}_{Z_\infty}/\partial\mathcal{O}_{Z_\infty}.$

\begin{thm}\label{thm:virob}
The obstruction appearing in Proposition~\ref{prop:virob} coincides,
up to the sign convention for the variational bracket, with the obstruction
computed from $\{\mu,\gamma\}$.
\end{thm}

\begin{proof}
Using the variational bracket, we have
$\{\mu,\gamma\}=
\int
\left(
\frac{\delta\mu}{\delta\theta}
\frac{\delta\gamma}{\delta\xi}
+
\frac{\delta\mu}{\delta\xi}
\frac{\delta\gamma}{\delta\theta}
\right).$
A direct computation gives
\[
\frac{\delta\mu}{\delta\theta}
=
-2\xi\theta^{(1)}-\xi^{(1)}\theta,
\qquad
\frac{\delta\mu}{\delta\xi}
=
\theta\theta^{(1)},
\]
and
therefore,
$
\{\mu,\gamma\}
=-4\int \xi\xi^{(1)}\theta\theta^{(1)}.
$
We identify $\xi^{(a)}\xi^{(b)}$ with $\partial^aL\wedge\partial^bL$, so that $\partial$ acts as $D_2$. Applying the symbol rule used in the proof of \cite[Theorem 4.6]{Fang26}, with the same ordering of the symbol variables and the second variable specialized to $-\lambda-D_2$, we obtain
\[
\int\xi\xi^{(1)}\theta\theta^{(1)}
\longmapsto
(2\lambda+D_2)(L\wedge\partial L)
=-(2\lambda+D_2)\delta(L).
\]
Hence $\{\mu,\gamma\}$ corresponds to $4(2\lambda+D_2)\delta(L)=4\delta([L_\lambda L])$, which is precisely the compatibility defect in Proposition~\ref{prop:virob}.
This is precisely the local functional representative of the failure of
the Virasoro $\lambda$-bracket and the cobracket $\delta$ to satisfy the
Lie conformal bialgebra compatibility condition. More explicitly,
Proposition~\ref{prop:virob} gives $L_\lambda\delta(L)-L_{-\lambda-D_2}\delta(L)
=
5\delta([L_\lambda L]),$
whereas compatibility would require the coefficient to be $1$. Hence the
obstruction is the remaining factor $4$, which is represented by the above
nonzero bracket $\{\mu,\gamma\}$.
\end{proof}

\section{Quantization of Lie conformal bialgebras via vertex bialgebras}\label{Sec:QT}

As a vertex algebraic analogue of ordinary bialgebras, Li introduced the notion of vertex bialgebras~\cite{Li07}, which may be regarded as a chiral counterpart of the bialgebraic structures underlying quantum groups. Subsequently, Linshaw and Malikov~\cite{LM19} constructed a vertex bialgebra quantizing the Khesin--Zakharevich Poisson--Lie group of pseudodifferential symbols. Their approach starts from the arc space $G_\infty$ of a Poisson--Lie group $G$ and considers the quantization of its coordinate vertex algebra $\mathcal O(G_\infty)$. 
In contrast, motivated by the quantized universal
enveloping (QUE) approach to ordinary Lie
bialgebras, we study $\hbar$-adic vertex bialgebra
deformations of universal enveloping vertex
algebras associated with Lie conformal algebras.
It would be interesting to investigate whether
these two approaches are related by a duality
analogous to that between quantized function
algebras and quantized universal enveloping
algebras in Poisson--Lie theory.

We now define what it means for an $\hbar$-adic
vertex bialgebra to quantize a Lie conformal
bialgebra.
 
A $\hbar$-adic vertex bialgebra quantizing a Lie conformal bialgebra
$(L,\delta)$ is a topologically free $\hbar$-adic vertex
bialgebra $A_\hbar$, together with an isomorphism of vertex
bialgebras $\Phi:A_\hbar/\hbar A_\hbar\simeq U(L)$,
such that the induced semiclassical cocommutator restricts
to $\delta$ on $L$. Thus, a $\hbar$-adic vertex bialgebra is
the \emph{chiral analogue of a quantized universal
enveloping algebra of an ordinary Lie bialgebra}.

\begin{thm}\label{thm:QUE}
Let $L$ be a Lie conformal algebra, and let $A:=U(L)$ carry its standard vertex bialgebra structure $(\Delta_0,\varepsilon_0)$. Let $(A_\hbar,\Delta_\hbar,\varepsilon_\hbar)$ be a topologically free $\hbar$-adic vertex bialgebra over $R:=\CC[[\hbar]]$, equipped with an isomorphism of vertex bialgebras $\Phi:A_\hbar/\hbar A_\hbar\xrightarrow{\sim}A.$
The vertex algebra structure and the counit are allowed to deform. For $a\in A$, choose a lift $\widetilde a\in A_\hbar$ satisfying $\Phi(\widetilde a\bmod\hbar)=a$, and define
\[
\delta_A(a):=(\Phi\otimes\Phi)\left(
\frac{\Delta_\hbar(\widetilde a)-\Delta_\hbar^{\mathrm{op}}(\widetilde a)}{\hbar}\bmod\hbar\right),
\qquad \Delta_\hbar^{\mathrm{op}}:=\tau\circ\Delta_\hbar,
\]
where $\tau(u\otimes v):=v\otimes u$. Then $\delta_A:A\to A\otimes A$ is well defined. If $\delta_A(L)\subseteq L\otimes L$, then $\delta:=\delta_A|_L$ makes $(L,[\cdot_\lambda\cdot],\delta)$ a Lie conformal bialgebra. Namely, $\delta$ is $\CC[\partial]$-linear, skew-symmetric, and satisfies the conformal $1$-cocycle condition and coJacobi identity
\[
\delta([x_\lambda y])=x_\lambda\delta(y)-y_{-\lambda-D_2}\delta(x),
\qquad
(\Id\otimes\delta)\delta-\tau_{12}(\Id\otimes\delta)\delta=(\delta\otimes\Id)\delta,
\]
where $D_2:=\partial\otimes\Id+\Id\otimes\partial$ and $\tau_{12}$ interchanges the first two tensor factors. 
Thus $A_\hbar$ is an $\hbar$-adic vertex
bialgebra quantization of
$(L,[\cdot_\lambda\cdot],\delta)$
in the sense defined above.

For the second order expansion, choose a vacuum preserving topological $R$-module isomorphism
$A_\hbar\simeq A[[\hbar]]$ inducing $\Phi$ modulo $\hbar$,
and write
$\Delta_\hbar=\Delta_0+\hbar\Delta_1+\hbar^2\Delta_2+O(\hbar^3)$,
where $\Delta_i\in\operatorname{Hom}_{\mathbb C}(A,A\otimes A)$.
Then the coefficient of $\hbar^2$ in coassociativity gives
\begin{align*}
&(\Delta_0\otimes\Id)\Delta_2+(\Delta_2\otimes\Id)\Delta_0
-(\Id\otimes\Delta_0)\Delta_2-(\Id\otimes\Delta_2)\Delta_0=(\Id\otimes\Delta_1)\Delta_1
-(\Delta_1\otimes\Id)\Delta_1.
\end{align*}
\end{thm}

\begin{proof}
Identify the special fiber with $A$ through $\Phi$, and choose a vacuum preserving topological $R$-module identification $A_\hbar\simeq A[[\hbar]]$. Write $Y_\hbar=Y_0+\hbar Y_1+\hbar^2Y_2+O(\hbar^3),\ 
\Delta_\hbar=\Delta_0+\hbar\Delta_1+\hbar^2\Delta_2+O(\hbar^3).$
Put $C_\hbar:=\Delta_\hbar-\Delta_\hbar^{\mathrm{op}}$. Since $\Delta_0$ is cocommutative, $C_\hbar(A_\hbar)\subseteq\hbar(A_\hbar\widehat\otimes_R A_\hbar).$
Topological freeness makes division by $\hbar$ well defined. If $\widetilde a'=\widetilde a+\hbar b$, then $\frac{C_\hbar(\widetilde a')-C_\hbar(\widetilde a)}{\hbar}\bmod\hbar
=\Delta_0(\overline b)-\Delta_0^{\mathrm{op}}(\overline b)=0.$
Hence $\delta_A$ is independent of the lift, and $\delta_A=\Delta_1-\tau\Delta_1,\ 
\tau\delta_A=-\delta_A,\ 
\delta_A(\mathbf1)=0.$
Compatibility with translation also gives $\delta_A\partial=D_2\delta_A$, where $D_2:=\partial\otimes\Id+\Id\otimes\partial$.

At first order, coassociativity gives $(\Delta_0\otimes\Id)\Delta_1+(\Delta_1\otimes\Id)\Delta_0
=(\Id\otimes\Delta_0)\Delta_1+(\Id\otimes\Delta_1)\Delta_0.$

For $x,y\in L$, choose lifts $\widetilde x,\widetilde y\in A_\hbar$. Since both $\Delta_\hbar$ and $\Delta_\hbar^{\mathrm{op}}$ are vertex algebra homomorphisms,
\begin{align*}
C_\hbar([\widetilde x_\lambda\widetilde y]_\hbar)
&=[\Delta_\hbar(\widetilde x){}_\lambda\Delta_\hbar(\widetilde y)]_\hbar
-[\Delta_\hbar^{\mathrm{op}}(\widetilde x){}_\lambda\Delta_\hbar^{\mathrm{op}}(\widetilde y)]_\hbar=[C_\hbar(\widetilde x){}_\lambda\Delta_\hbar(\widetilde y)]_\hbar
+[\Delta_\hbar^{\mathrm{op}}(\widetilde x){}_\lambda C_\hbar(\widetilde y)]_\hbar.
\end{align*}
Dividing by $\hbar$ and reducing modulo $\hbar$ gives $\delta([x_\lambda y])
=[\delta(x){}_\lambda\Delta_0(y)]
+[\Delta_0(x){}_\lambda\delta(y)].$
Since $\Delta_0(x)=x\otimes\mathbf1+\mathbf1\otimes x$, the second term is $x_\lambda\delta(y)$, where $x_\lambda(u\otimes v):=[x_\lambda u]\otimes v+u\otimes[x_\lambda v].$
Skew-symmetry in the tensor product vertex algebra therefore gives $\delta([x_\lambda y])
=x_\lambda\delta(y)-y_{-\lambda-D_2}\delta(x).$ Moreover, as in the Lie bialgebra case, the coassociativity
of $\Delta_\hbar$ implies the coJacobi identity $(\Id\otimes\delta)\delta
-\tau_{12}(\Id\otimes\delta)\delta
=(\delta\otimes\Id)\delta.$

At second order, comparing the coefficients of $\hbar^2$ in coassociativity yields
\begin{align*}
&(\Delta_0\otimes\Id)\Delta_2+(\Delta_2\otimes\Id)\Delta_0
-(\Id\otimes\Delta_0)\Delta_2-(\Id\otimes\Delta_2)\Delta_0\\
&\qquad=(\Id\otimes\Delta_1)\Delta_1-(\Delta_1\otimes\Id)\Delta_1.
\end{align*}
To extract the coJacobi identity, use the cocommutator identity implied by coassociativity:
\[
(\Id\otimes C_\hbar)C_\hbar-\tau_{12}(\Id\otimes C_\hbar)C_\hbar
=(C_\hbar\otimes\Id)C_\hbar.
\]
Since $C_\hbar=\hbar\delta_A+O(\hbar^2)$, its $\hbar^2$ coefficient is $(\Id\otimes\delta_A)\delta_A-\tau_{12}(\Id\otimes\delta_A)\delta_A
=(\delta_A\otimes\Id)\delta_A.$
Restricting to $L$ gives the co-Jacobi identity for $\delta$, completing the proof.
\end{proof}

\begin{rmk}\label{rmk:QUE0}
Indeed, the condition $\delta_A(L)\subseteq L\otimes L$ is automatic. Keeping the above notation, the coefficient of $\hbar$ in $(\Delta_\hbar\otimes\Id)\Delta_\hbar=(\Id\otimes\Delta_\hbar)\Delta_\hbar$ gives
\[
(\Delta_0\otimes\Id)\Delta_1+(\Delta_1\otimes\Id)\Delta_0
=(\Id\otimes\Delta_0)\Delta_1+(\Id\otimes\Delta_1)\Delta_0.
\]
Moreover, the permutation identity $(\Id\otimes\Delta_\hbar)\tau\Delta_\hbar=\tau_{12}\tau_{23}(\Delta_\hbar\otimes\Id)\Delta_\hbar$ and coassociativity imply
\[
(\Id\otimes\Delta_\hbar)C_\hbar=(C_\hbar\otimes\Id)\Delta_\hbar+\tau_{12}(\Id\otimes C_\hbar)\Delta_\hbar,
\]
where $C_\hbar:=\Delta_\hbar-\tau\Delta_\hbar$. Since $C_\hbar=\hbar(\Delta_1-\tau\Delta_1)+O(\hbar^2)=\hbar\delta_A+O(\hbar^2)$, comparing the coefficients of $\hbar$ gives
\begin{equation}\label{eq:QUE0}
(\Id\otimes\Delta_0)\delta_A=(\delta_A\otimes\Id)\Delta_0+\tau_{12}(\Id\otimes\delta_A)\Delta_0.
\end{equation}
For $x\in L$, since $\Delta_0(x)=x\otimes\vac+\vac\otimes x$ and $\delta_A(\vac)=0$, identity~\eqref{eq:QUE0} gives $(\Id\otimes\Delta_0)\delta_A(x)=\delta_A(x)\otimes\vac+\tau_{12}(\vac\otimes\delta_A(x)).$
Now, write $\delta_A(x)=\sum_i a_i\otimes b_i\in U(L)\otimes U(L)$, where the $a_i$ are linearly independent. Then $\sum_i a_i\otimes\bigl(\Delta_0(b_i)-b_i\otimes\vac-\vac\otimes b_i\bigr)=0.$
By linear independence, $\Delta_0(b_i)=b_i\otimes\vac+\vac\otimes b_i$ for every $i$. Thus each $b_i$ is primitive and hence belongs to $L$, since the primitive elements of $U(L)$ are precisely those of $L$. Therefore $\delta_A(x)\in U(L)\otimes L$. Together with $\tau\delta_A=-\delta_A$, this implies $\delta_A(x)$ in $L\otimes L.$
\end{rmk}

\begin{rmk}\label{rmk:QUE1}
We may also deform the vertex algebra structure. Under the above identification, write $Y_\hbar(a,z)=\sum_{n\geq0}\hbar^nY_n(a,z)$. Since $\Delta_\hbar$ is a vertex algebra homomorphism, we have $\Delta_\hbar(Y_\hbar(a,z)b)=Y^\otimes_\hbar(\Delta_\hbar(a),z)\Delta_\hbar(b),$
where $Y_\hbar^\otimes$ is the tensor product vertex operator determined by $Y_\hbar^\otimes(a\otimes b,z)(c\otimes d)
=Y_\hbar(a,z)c\otimes Y_\hbar(b,z)d.$
Thus $Y_\hbar^\otimes$ is completely determined by $Y_\hbar$ and requires no independent deformation data. This is related to the first order deformation theory of vertex algebras studied in~\cite{KQ24}, with the additional requirement that the deformed vertex algebra structure be compatible with the coproduct $\Delta_\hbar$.
\end{rmk}

\begin{thm}\label{thm:CQ}
Let $L=\operatorname{Cur}\mathfrak{sl}_2$, and let
$A=U(L)$ be the universal level zero affine vertex algebra
with vacuum $\mathbf1$, translation operator $\partial$.
Fix an odd polynomial $p(t)\in\mathbb C[t]$, and set
$$
q_p=p(\partial)h,\qquad
K_{\hbar,p}
=\exp\!\left(\frac{\hbar}{2}q_p\right)
=\sum_{m\geq0}\frac{\hbar^m}{2^m m!}q_p^m,
$$
where the powers are taken in the commutative differential
vertex subalgebra generated by $h$, with $q_p^0=\mathbf1$.

Equip $A_\hbar=A[[\hbar]]$ with its coefficientwise vertex
operations, and use the completed tensor product $A_\hbar\widehat\otimes A_\hbar
=(A\otimes_{\mathbb C}A)[[\hbar]].$
All vertex algebra identities are understood $\hbar$-adically,
equivalently modulo every positive power of $\hbar$.

Then the assignments
\begin{align*}
\Delta_{\hbar,p}(h)
&=h\otimes\mathbf1+\mathbf1\otimes h,\\
\Delta_{\hbar,p}(e)
&=e\otimes K_{\hbar,p}
  +K_{\hbar,p}^{-1}\otimes e,\\
\Delta_{\hbar,p}(f)
&=f\otimes K_{\hbar,p}^{-1}
  +K_{\hbar,p}\otimes f
\end{align*}
extend uniquely to a continuous $\mathbb C[[\hbar]]$-linear
vertex algebra homomorphism
$
\Delta_{\hbar,p}:A_\hbar
\longrightarrow A_\hbar\widehat\otimes A_\hbar.
$
Together with the standard augmentation $\varepsilon:A_\hbar\longrightarrow\mathbb C[[\hbar]],
\ 
\varepsilon(\mathbf1)=1,\ 
\varepsilon(e)=\varepsilon(h)=\varepsilon(f)=0,$
this makes $A_\hbar$ a topologically free $\hbar$-adic
vertex bialgebra. In particular,
\begin{gather*}
(\Delta_{\hbar,p}\otimes\Id)\Delta_{\hbar,p}
=(\Id\otimes\Delta_{\hbar,p})\Delta_{\hbar,p},\ 
(\varepsilon\otimes\Id)\Delta_{\hbar,p}
=(\Id\otimes\varepsilon)\Delta_{\hbar,p}
=\Id.
\end{gather*}

Moreover, $\Delta_{\hbar,p}\equiv\Delta_0\pmod{\hbar}$,
where $\Delta_0$ is the standard primitive coproduct.
Writing $\tau(a\otimes b)=b\otimes a$ and
$a\wedge b=a\otimes b-b\otimes a$, the classical limit
$$
\delta_p(a)
=\frac{\Delta_{\hbar,p}(a)-\tau\Delta_{\hbar,p}(a)}
       {\hbar}\bmod\hbar,
\qquad a\in L,
$$
is the Lie conformal cobracket determined by $\delta_p(h)=0,\ 
\delta_p(e)=e\wedge p(\partial)h,\ 
\delta_p(f)=-f\wedge p(\partial)h.$
Thus $(A_\hbar,\Delta_{\hbar,p},\varepsilon)$ is an
all orders quantization of $(L,\delta_p)$ with the
underlying vertex algebra fixed. For $p\neq0$, its
coproduct is noncocommutative. the coboundary cobracket associated with the conformal $r$-matrix $r_p$ in Example~\ref{exa:RBJCa}.
\end{thm}

\begin{proof}
We give the proof for $p(t)=t$. For a general odd polynomial $p$, the same argument applies with $q=p(\partial)h$, using $p(-\partial)=-p(\partial)$. Set $q:=\partial h,\qquad D:=(h\otimes q)_{(0)},\ \Phi_\hbar:=\exp\left(\frac{\hbar}{4}D\right),\  \Delta_\hbar:=\Phi_\hbar\circ\Delta_0.$
Since every zero mode is a derivation, $\Phi_\hbar$ is a continuous vertex algebra automorphism of $(A\otimes A)[[\hbar]]$. Hence $\Delta_\hbar$ is a continuous vertex algebra homomorphism. A direct computation gives
\[
D(e\otimes\vac)=2e\otimes q,\quad D(\vac\otimes e)=-2q\otimes e,\quad D(f\otimes\vac)=-2f\otimes q,\quad D(\vac\otimes f)=2q\otimes f,
\]
and $D(h\otimes\vac)=D(\vac\otimes h)=0$. Since $[q_\lambda q]=0$, the vertex subalgebra generated by $q$ is commutative; we write $q^m:=:q\cdots q:$ and set $K:=\exp\left(\frac{\hbar}{2}q\right)=\sum_{m\geq0}\frac{1}{m!}\left(\frac{\hbar}{2}\right)^m q^m.$
It follows that
\[
\Delta_\hbar(h)=h\otimes\vac+\vac\otimes h,\qquad \Delta_\hbar(e)=e\otimes K+K^{-1}\otimes e,\qquad \Delta_\hbar(f)=f\otimes K^{-1}+K\otimes f.
\]
Moreover, $\Delta_\hbar(q)=q\otimes\vac+\vac\otimes q.$
Since $:(a\otimes\vac)(\vac\otimes b):=a\otimes b$ and the subalgebra generated by $q\otimes\vac$ and $\vac\otimes q$ is commutative, the ordinary binomial formula gives $\Delta_\hbar(K)=K\otimes K,\  \Delta_\hbar(K^{-1})=K^{-1}\otimes K^{-1}.$
Hence coassociativity follows by checking the generators $e,h,f$; for example,
\[
(\Delta_\hbar\otimes\Id)\Delta_\hbar(f)=f\otimes K^{-1}\otimes K^{-1}+K\otimes f\otimes K^{-1}+K\otimes K\otimes f=(\Id\otimes\Delta_\hbar)\Delta_\hbar(f),
\]
and the cases of $e$ and $h$ are similar.

Finally, let $\epsilon:A[[\hbar]]\to\CC[[\hbar]]$ be the standard augmentation, with $\epsilon(\vac)=1$ and $\epsilon(e)=\epsilon(h)=\epsilon(f)=0$. Since $\epsilon(q)=0$, we have $\epsilon(K^{\pm1})=1$. Thus, for example,
\[
(\epsilon\otimes\Id)\Delta_\hbar(f)=1\otimes f\cong f,\qquad (\Id\otimes\epsilon)\Delta_\hbar(f)=f\otimes1\cong f,
\]
and similarly for $e$ and $h$. Therefore
\[
(\Delta_\hbar\otimes\Id)\Delta_\hbar=(\Id\otimes\Delta_\hbar)\Delta_\hbar,\qquad (\epsilon\otimes\Id)\Delta_\hbar=\Id=(\Id\otimes\epsilon)\Delta_\hbar.
\]
Finally, $A[[\hbar]]$ is topologically free over $\CC[[\hbar]]$, and $K^{\pm1}=1\pm\frac{\hbar}{2}q+O(\hbar^2),$
so $\Delta_\hbar\equiv\Delta_0\pmod{\hbar}$. Moreover, $\frac{\Delta_\hbar-\Delta_\hbar^{\op}}{\hbar}\bmod\hbar$
vanishes on $h$ and gives $\delta(e)=e\wedge q,\  \delta(f)=-f\wedge q.$
Thus $\Delta_\hbar$ is an all orders quantization of the corresponding Lie conformal bialgebra.
\end{proof}
It would be interesting to study tensor products of Wakimoto modules for $V^0(\fsl_2)$ with respect to this coproduct.

\begin{rmk}\label{rmk:GCQ}
The zero mode construction in the proof of Theorem~\ref{thm:CQ} works
uniformly for every odd polynomial $p$. Indeed, set $q_p=p(\partial)h,\  D_p=(h\otimes q_p)_{(0)}.$
Using $p(-\partial)=-p(\partial)$, we obtain
\[
D_p(a\otimes1)=[h,a]\otimes q_p,\qquad
D_p(1\otimes a)=-q_p\otimes[h,a],
\quad a\in\mathfrak{sl}_2.
\]
For example, $D_p(1\otimes e)
=\sum_{j\geq0}\frac{\partial^j h}{j!}\otimes(q_p)_{(j)}e
=2p(-\partial)h\otimes e
=-2q_p\otimes e.$
Since $D_p$ is a vertex algebra derivation, $\Phi_{\hbar,p}
=\exp\!\left(\frac{\hbar}{4}D_p\right)$
defines a continuous vertex algebra automorphism of
$(A\otimes A)[[\hbar]]$. The remainder of the proof proceeds exactly as in
Theorem~\ref{thm:CQ}.
\end{rmk}

\begin{rmk}\label{rmk:BFfusion}
Motivated by the holomorphic-topological BF theory
in~\cite[Eq.~(5.11)]{Zeng26}, the Cartan type
cobracket $\delta_p$ suggests the candidate
interaction functional $S_p^{\mathrm{BF}}
=\int_M (A^eB_e-A^fB_f)\,
p(\partial_w)B_h\wedge dw,$
where $M=\mathbb R^2\times\mathbb C$, $w$ is the
holomorphic coordinate, and $A^a,B_a$ are the
degree one BF superfields associated with the
basis $a=e,h,f$ and its dual.
Here $p(\partial_w)$ acts on $B_h$, and the
integral selects the top form component,
formally modulo total derivatives.
A precise identification of the corresponding
BV derived operation with $\delta_p$ remains
to be established.
\end{rmk}




\begin{thebibliography}{99}
\raggedright

\bibitem{AN25}
R.~Abedin and W.~Niu,
\emph{Double Yangian, factorization, and qKZ-equation for cotangent Lie algebras},
arXiv:2512.16996 [math.QA], 2025.

\bibitem{AS05}
A.~Alekseev and T.~Strobl,
\emph{Current algebras and differential geometry},
J. High Energy Phys. \textbf{2005} (2005), no.~03, 035.

\bibitem{ACKL17}
T.~Arakawa, T.~Creutzig, K.~Kawasetsu, and A.~R.~Linshaw,
\emph{Orbifolds and cosets of minimal $\mathcal{W}$-algebras},
Commun. Math. Phys. \textbf{355} (2017), no.~1, 339--372.

\bibitem{AM}
T.~Arakawa and A.~Moreau,
\emph{Arc spaces and vertex algebras},
monograph, online version dated October 25, 2024.

\bibitem{BDSK20}
B.~Bakalov, A.~De Sole, and V.~G.~Kac,
\emph{Computation of cohomology of Lie conformal and Poisson vertex algebras},
Selecta Math. (N.S.) \textbf{26} (2020), no.~4, Paper No.~50, 51~pp.

\bibitem{BKV99}
B.~Bakalov, V.~G.~Kac, and A.~A.~Voronov,
\emph{Cohomology of conformal algebras},
Commun. Math. Phys. \textbf{200} (1999), no.~3, 561--598.

\bibitem{BDSK09}
A.~Barakat, A.~De Sole, and V.~G.~Kac,
\emph{Poisson vertex algebras in the theory of Hamiltonian equations},
Jpn. J. Math. \textbf{4} (2009), no.~2, 141--252.

\bibitem{BG22}
C.~Boyallian and J.~Guzm\'an,
\emph{Formal vertex laws associated to Lie conformal algebras},
J. Math. Phys. \textbf{63} (2022), no.~7, 071701.

\bibitem{BN26}
D.~Butson and S.~Nair,
\emph{On the deformation theory of chiral quantizations},
arXiv:2606.27341 [math.QA], 2026.

\bibitem{CL22}
T.~Creutzig and A.~R.~Linshaw,
\emph{Trialities of $\mathcal{W}$-algebras},
Camb. J. Math. \textbf{10} (2022), no.~1, 69--194.

\bibitem{DSK05}
A.~De Sole and V.~G.~Kac,
\emph{Freely generated vertex algebras and non-linear Lie conformal algebras},
Commun. Math. Phys. \textbf{254} (2005), no.~3, 659--694.

\bibitem{DSK06}
A.~De Sole and V.~G.~Kac,
\emph{Finite vs.\ affine $\mathcal{W}$-algebras},
Jpn. J. Math. \textbf{1} (2006), no.~1, 137--261.

\bibitem{DSK09}
A.~De Sole and V.~G.~Kac,
\emph{Lie conformal algebra cohomology and the variational complex},
Commun. Math. Phys. \textbf{292} (2009), no.~3, 667--719.

\bibitem{Dr83}
V.~G.~Drinfeld,
\emph{Hamiltonian structures on Lie groups, Lie bialgebras and the
geometric meaning of the classical Yang--Baxter equations},
Soviet Math. Dokl. \textbf{27} (1983), no.~1, 68--71.

\bibitem{Dr87}
V.~G.~Drinfeld,
\emph{Quantum groups},
in \emph{Proceedings of the International Congress of Mathematicians}
(Berkeley, 1986), Vol.~1,
American Mathematical Society, Providence, RI, 1987,
pp.~798--820.

\bibitem{EK96}
P.~Etingof and D.~Kazhdan,
\emph{Quantization of Lie bialgebras, I},
Selecta Math. (N.S.) \textbf{2} (1996), no.~1, 1--41.

\bibitem{EK98}
P.~Etingof and D.~Kazhdan,
\emph{Quantization of Lie bialgebras, II},
Selecta Math. (N.S.) \textbf{4} (1998), no.~2, 213--231.

\bibitem{Fang25}
W.~Fang,
\emph{Generalized AKS scheme of integrability via vertex algebra},
Int. Math. Res. Not. IMRN \textbf{2025} (2025), no.~11, rnaf143.

\bibitem{FangPhD}
W.~Fang,
\emph{On the Chiralization of Poisson--Lie Theory with Applications
in Integrable Systems},
Ph.D. thesis, Kyoto University, 2026.

\bibitem{Fang26}
W.~Fang,
\emph{Shifted symplectic structures and Poisson vertex algebra},
arXiv:2601.17840v2 [math.QA], 2026.

\bibitem{Fangpre}
W.~Fang,
\emph{Geometric theory of chiralization of Lie bialgebras},
in preparation.

\bibitem{GK23}
V.~Gubarev and R.~Kozlov,
\emph{Conformal Yang--Baxter equation on
$\mathrm{Cur}(\mathfrak{sl}_2(\mathbb{C}))$},
J. Math. Phys. \textbf{64} (2023), no.~1, 011704.

\bibitem{GKRB23}
V.~Gubarev and R.~Kozlov,
\emph{Rota--Baxter operators on
$\mathrm{Cur}(\mathfrak{sl}_2(\mathbb{C}))$},
Int. Electron. J. Algebra \textbf{33} (2023), 247--269.

\bibitem{HLX22}
J.~Han, H.~Li, and Y.~Xiao,
\emph{Cocommutative vertex bialgebras},
J. Algebra \textbf{598} (2022), 536--569.

\bibitem{HB20}
Y.~Hong and C.~Bai,
\emph{Conformal classical Yang--Baxter equation, $S$-equation and
$\mathcal{O}$-operators},
Lett. Math. Phys. \textbf{110} (2020), no.~5, 885--909.

\bibitem{IK13}
N.~Ikeda and K.~Koizumi,
\emph{Current algebras and QP-manifolds},
Int. J. Geom. Methods Mod. Phys. \textbf{10} (2013), no.~6, 1350024.

\bibitem{IX13}
N.~Ikeda and X.~Xu,
\emph{Canonical functions, differential graded symplectic pairs
in supergeometry, and Alexandrov--Kontsevich--Schwartz--Zaboronsky
sigma models with boundaries},
J. Math. Phys. \textbf{55} (2014), no.~11, 113505.

\bibitem{IX13v}
N.~Ikeda and X.~Xu,
\emph{Current algebras from QP-manifolds in general dimensions},
J. Math. Phys. \textbf{66} (2025), no.~1, 012301.

\bibitem{Kac98}
V.~G.~Kac,
\emph{Vertex Algebras for Beginners},
2nd ed., University Lecture Series, vol.~10,
American Mathematical Society, Providence, RI, 1998.

\bibitem{K15}
V.~G.~Kac,
\emph{Introduction to vertex algebras, Poisson vertex algebras,
and integrable Hamiltonian PDE},
in \emph{Perspectives in Lie Theory},
Springer INdAM Series, vol.~19,
Springer, Cham, 2017, pp.~3--72.

\bibitem{KS11}
Y.~Kosmann-Schwarzbach,
\emph{Nijenhuis structures on Courant algebroids},
Bull. Braz. Math. Soc. (N.S.) \textbf{42} (2011), no.~4, 625--649.

\bibitem{KQ24}
V.~Kovalchuk and F.~Qi,
\emph{First-order deformations of freely generated vertex algebras},
arXiv:2408.16309 [math.QA], 2024.

\bibitem{Li04}
H.~Li,
\emph{Vertex algebras and vertex Poisson algebras},
Commun. Contemp. Math. \textbf{6} (2004), no.~1, 61--110.

\bibitem{Li07}
H.~Li,
\emph{A smash product construction of nonlocal vertex algebras},
Commun. Contemp. Math. \textbf{9} (2007), no.~5, 605--637.

\bibitem{Lib08}
J.~I.~Liberati,
\emph{On conformal bialgebras},
J. Algebra \textbf{319} (2008), no.~6, 2295--2318.

\bibitem{LM19}
A.~Linshaw and F.~Malikov,
\emph{One example of a chiral Lie group},
to appear in RIMS K\^oky\^uroku Bessatsu.
arXiv:1902.07414, 2019.

\bibitem{LWX97}
Z.-J.~Liu, A.~Weinstein, and P.~Xu,
\emph{Manin triples for Lie bialgebroids},
J. Differential Geom. \textbf{45} (1997), no.~3, 547--574.

\bibitem{Primc99}
M.~Primc,
\emph{Vertex algebras generated by Lie algebras},
J. Pure Appl. Algebra \textbf{135} (1999), no.~3, 253--293.

\bibitem{Roy01}
D.~Roytenberg,
\emph{Courant algebroids, derived brackets and even symplectic
supermanifolds},
Ph.D. thesis, University of California, Berkeley, 1999.
arXiv:math/9910078 [math.DG].

\bibitem{Vor01}
T.~Voronov,
\emph{Graded manifolds and Drinfeld doubles for Lie bialgebroids},
in \emph{Quantization, Poisson Brackets and Beyond}
(Manchester, 2001),
Contemp. Math., vol.~315,
American Mathematical Society, Providence, RI, 2002,
pp.~131--168.

\bibitem{Yam07}
G.~Yamskulna,
\emph{Vertex Poisson algebras associated with Courant algebroids},
Comm. Algebra \textbf{35} (2007), no.~12, 4137--4162.

\bibitem{Zeng26}
K.~Zeng,
\emph{Factorization algebras and quantum groups from generalized Poisson sigma models},
arXiv:2607.09486, 2026.

\end{thebibliography}
\end{document}